\documentclass[10pt,openbib]{article}  
\usepackage[square,numbers,sort&compress]{natbib} 

\newlength{\bibitemsep}
\newlength{\bibparskip}
\let\oldthebibliography\thebibliography
\renewcommand\thebibliography[1]{%
	\oldthebibliography{#1}%
	\setlength{\parskip}{2.1\bibitemsep}%
	\setlength{\itemsep}{2.0\bibparskip}%
}

\usepackage[export]{splitbib}  
 
\begin{category}[]{\normalsize \mbox{} \vspace{-2mm} \mbox{} \\ -- \hspace{0mm} e-value \ Theory and Epistemology}
\SBentries{Borges2006,Borges2007,Cabras2015,Esteves2016,Esteves2019,Esteves2023,Fossaluza2017,Izbicki2015,Izbicki2022,KelterStern2020,Kelter2021b,Lauretto2005a,Lauretto2007,Madruga2001,Madruga2003,PereiraStern1999,Pereira2008,Pereira2011,PereiraStern2020,Pinto2012,Ruli2020,Ruli2020b,Silva2015,Silva2018,Stern2003,Stern2004,Stern2007a,Stern2007b,Stern2008a,Stern2008b,Stern2011a,Stern2011b,Stern2014a,Stern2014b,SternPereira2014,Stern2015,Stern2017a,Stern2017b,Stern2018a,Stern2018b,Stern2018c,Stern2020,Stern2022,Stern2023,Ventura2013,Ventura2014,Ventura2016}
\end{category} 
   
\begin{category}[]{\normalsize  \mbox{} \vspace{-2mm} \mbox{} \\ --  \hspace{0mm} e-value \ Use in Empirical Science}
\SBentries{Ainsbury2013,Andrade2015,Andrade2014,Assane2018,Assane2019,Barahona2016,Bernardini2011,Bernardo2012,Brentani2011,Cantinha2017,Camargo2012,Cerezetti2012,Chakrabarty2017,Chaiboonsri2018,Chen2015,DCunha2016,Diniz2011,Diniz2012a,Diniz2012b,Diniz2022,Hubert2009,Garcia2016,Hubert2018,Hubert2018b,Hubert2019,Irony2002,Izbicki2012,Johnson2009,Kelter2020,Kelter2021a,Kostrzewski2012,Lauretto2003,Lauretto2005b,Lauretto2009,Lima2014,Loschi2007,Loschi2012,Maranhao2012,Mathis2008,Montoya2001,Maranhao201,Nakano2006,Oliveira2018,Patriota2013,Patriota2017,PereiraStern1999b,Ranzato2018,Rifo2009,Rifo2012,Rincon2010,Rodrigues2006,Santos2020,Seixas2008,Shavitt2017,Sikov2019,Spektor2019,SternZacks2002,Takada2020,Vikas2016,Vosseler2016,Wittenburg2016}
\end{category}   

\begin{category}[]{\normalsize  \mbox{} \vspace{-2mm} \mbox{} \\ -- \hspace{0mm} General References in Logic and Statistics}
\SBentries{Aczel1975,Amari1987,Amari2007,Amari2016,Andrews2001,Angus1994,Attneave1959,Barlow1981,Barlow1993,Basu1988,Benthem2000,Berger1988,Bernardo2005,Beziau2012,Birgin2014,Birnbaum1961,Blanche1966,Bonassi2009,Borges1979,Box1973,Bracewell1986,Brooks2011,Carnielli2008a,Carnielli2008c,Carnielli2010,Chen2000,Cherkassky1998,Chernoff1954,Cohen1977,Costa2001,Darwiche1992,Darwiche1993,DeGroot1970,DeGroot2011,Drton2009,Dubois1982,Dubois2004,Dubois2007,Dubois2012,Dugdale1996,Dugundji1966,Eaton1989,Ebbinghaus1999,Epstein1974,Epstein1993,Epstein2008,Evans1997,Fang1997,Finetti1957,Finetti1975,Flanders1974,Fong2019,Frieden2004,Gaifman1964,Gaifman1982,Gelman2004,Gilks1996,Gottenwald2001,Guy1994,Hailperin1996,Hailperin2011,Hammersley1964,Hedman2004,Hocking1985,Howson2009,Hutter2005,Hutter2008,Inhasz2010,Jeffreys1961,Kadane2011,Kadane2016,Kaplan1987,Kapur1989,Kapur1992,Karawatzki2005,Kaufmann1977,Kelley1987a,Kelley1987b,Klir1988,Klir1995,Lawvere1997,Lindley1972,Lindley1991,Lipschutz1965,Luenberger2008,Malinowski1993,Minoux1986,Pawitan2001,PereiraNS2000,Pereira2a,PereiraStern2008b,Pereira2b,Pflug1996,Piskunov1969,Priest2001,Rine1984,Ripley1987,Rissanen1989,Rodrigues2021,Romero1991,Royall1997,Ruhl2008,Ruhl2015,Russell2015,Sakamoto1986,Salicone2007,Schwartzman2008,Searcoid2006,Self1987,Shackle1968,Shackle1969,Shier1991,Sichelman2021,Skilling2004,Skilling2006,Spall2003,Spivak2014,Springer1979,Stern1992,SternColla2009,Stern2009,SternVavasis1994,Thulin2014,Vieland2019,Voigt1990,Wah2008,Wallace1968,Wallace1999,Wallace2005,Wechsler2008,West1997,Williams2021,Williamson1989,Zeleny1982,Zellner1971} 
\end{category}   

\begin{category}[]{\normalsize  \mbox{} \vspace{-2mm} \mbox{} \\ -- \hspace{0mm} General References in Law and Philosophy}
\SBentries{Barak2012,Barker2001,Bayes1764,Bueno2016,Campilongo2020,Cauble2021,Cobde1998,Dumitriu1977,Engle2012,Foerster2001,Foerster2003,Finetti2006,Gallais1982,Garson2006,Gaskins1992,Good1983,Hulsroj2013,Kades1997,Kokott1998,Luhmann1985,Luhmann1989,Maturana1980,Muchmore2016,Murray2019,Pigliucci2014,Rasch2000,Schuck1992,Teubner1988,Wittgenstein1921,Zeleny1979}
\end{category}

\usepackage{geometry}
\usepackage{multicol} 
\usepackage{mathptmx}      

\usepackage{graphicx}

\usepackage{amsmath} 
\usepackage{mathrsfs}  
\usepackage{amssymb} 
 
\def\overli{\overline} 
 
\usepackage{amsfonts}
\usepackage{fancyhdr}

\def\argmin \mathop{\rm argmin}
\def\Re{I\kern -0.37 em R}
\def\Na{I\kern -0.37 em N}
\def\Qe{I\kern -0.37 em Q}

\def\g{\,\vert \,}

\def\ev{\,\mbox{ev}\,} 
\def\evb{\overline{\ev}} 
\def\sev{\mbox{sev}} 
\def\sevb{\overline{\sev}} 
 
\def\zero{\mbox{\bf 0}}

\def\Chi2{\mbox{Q}} 
\def\pprod{\prod\nolimits}

\def\mmax{\max\nolimits} 
\def\bbigvee{\bigvee\nolimits} 
\def\bbigwedge{\bigwedge\nolimits} 
\def\g{\,\vert\,}

\def\Pr{\mbox{Pr}}

\def\Tb{\overli{T}}

\def\Hb{\overli{H}}
\def\Wb{\overli{W}}

\def\zero{\mbox{\bf 0}}

\def\Re{{R}}

\begin{document}

   
\title{  
The e-value and the Full Bayesian Significance Test: \\ 
Logical Properties and Philosophical Consequences.
	   
         
\vspace{3mm} \mbox{} }


\author{
{Julio Michael Stern\thanks{ 
\texttt{jstern@ime.usp.br\,, www.ime.usp.br/\~{}jmstern/}\,, 
ORCID: 0000-0003-2720-3871, 
ResearcherID: C-1128-2013,  
ID-Lattes: 9582404119292455, 
Google Scholar: 7spXyx8AAAAJ; \ 
IME-USP - 
Institute of Mathematics and Statistics of the University of Sao Paulo. 
Rua do Mat\~{a}o, 1010, 05508-900, S\~{a}o Paulo, Brazil.  
}
\hspace{4mm} \mbox{} } 
\and 
{Carlos Alberto de Bragan\c{c}a Pereira\thanks{ 
\texttt{cpereira@ime.usp.br}, 
ORCID: 0000-0001-6315-7484, 
ResearcherID: O-5022-2015; IME-USP. 
}
\hspace{0mm} \mbox{} } 
\and 
{Marcelo de Souza Lauretto\thanks{
\texttt{marcelolaureto@usp.br}, 
ORCID: 0000-0001-5507-2368, 
ResearcherID: J-1964-2012;  
EACH-USP - School of Arts, Sciences and Humanities of the University of Sao Paulo. 
Rua Arlindo B\'{e}ttio, 1000, 03828-000, S\~{a}o Paulo, Brazil. 
}
\hspace{16mm} \mbox{} } 
\and
{Lu\'{\i}s Gustavo Esteves\thanks{ 
\texttt{lesteves@ime.usp.br}, 
ID-Lattes: 8881931568548501; IME-USP. 
}
\hspace{0mm} \mbox{} } 
\and 
{Rafael Izbicki\thanks{ 
\texttt{rizbicki@ufscar.br}, 
ORCID: 0000-0003-0379-9690, 
ResearcherID: J-1160-2014; 
DS-UFScar - Department of Statistics of the Federal University of Sao Carlos. 
Rodovia Washington Luis, km 235, 13565-905, S\~{a}o Carlos, Brazil. 
}
\hspace{43mm} \mbox{} } 
\and
{Rafael Bassi Stern\thanks{ 
\texttt{rbstern@ufscar.br}, 
ORCID: 0000-0002-4323-4515, 
ResearcherID: I-4745-2014;  
DS-UFSCar. 
}
\hspace{0mm} \mbox{} } 
\and
{Marcio Alves Diniz\thanks{ 
\texttt{marcio.alves.diniz@gmail.com}, 
ORCID: 0000-0002-8239-4263, 
ID-Lattes: 8948404469003829;  
DS-UFSCar. 
}
\hspace{22mm} \mbox{} } 
\and
{Wagner de Souza Borges\thanks{ 
\texttt{wborges@ime.usp.br}, 
ORCID: 0000-0001-5520-5861,  
ID-Lattes: 0975898314511610;   
IME-USP. 
}
\hspace{0mm} \mbox{} } 
%
}

\date{} 

\maketitle 

\mbox{} \vspace{-2mm} \mbox{} 
\begin{center} 
\textit{Dedicated to Walter Alexandre Carnielli, for his 70th birthday}  
\end{center}

\mbox{} \\ 
\vspace{-3mm} \mbox{}

\begin{quotation}
\textbf{Abstract:} This article gives a conceptual review of the e-value, $\ev(H\g X)$ -- the epistemic value of hypothesis $H$ given observations $X$. 
This statistical significance measure was developed in order to allow logically coherent and consistent tests of hypotheses, including sharp or precise hypotheses, via the Full Bayesian Significance Test (FBST). 
Arguments of analysis allow a full characterization of this statistical test by its logical or compositional properties, showing a mutual complementarity between results of mathematical statistics and the logical desiderata lying at the foundations of this theory.  


\end{quotation} 

\mbox{} \vspace{-4mm} \mbox{} 

\section{Introduction} 
\label{intro}

The e-value, $\mbox{ev}(H|X)$ -- also named the \textit{epistemic-value} of hypothesis H given observations X, or the evidence-value of observations X in favor (or in support) of hypothesis H -- is a Bayesian statistical significance measure introduced in 1999 by Carlos Alberto de Bragan\c{c}a Pereira and Julio Michael Stern, together with the FBST -- the Full Bayesian Significance Test, see \cite{PereiraStern1999}. 
The definitions of e-value and the FBST were further refined and generalized by subsequent works of several researchers at the University of S\~{a}o Paulo (USP) and the Federal University of S\~{a}o Carlos (UFSCar), in Brazil, including Wagner Borges, Lu\'{\i}s Gustavo Esteves, Rafael Izbicki, Regina Madruga, Rafael Bassi Stern, and Sergio Wechsler, see \cite{Borges2007,Esteves2016,Esteves2019,Madruga2003,Pereira2008}. 

The e-value was specially designed to assess the \textit{logical truth value} (a.k.a. \textit{statistical significance})  of \textit{sharp} (a.k.a. \textit{precise}) hypotheses in the context of Bayesian statistics. 
The e-value has desirable asymptotic, geometrical (i.e.,\,invariance), and logical (i.e.,\,compositional) properties that allow consistent and coherent evaluation and testing of sharp statistical hypotheses. 
Furthermore, in applied modeling, the FBST offers an easy-to-implement and powerful statistical test that is fully compliant with Bayesian principles of good inference, like the likelihood principle, see \cite{Berger1988,Gelman2004,Pawitan2001,Wechsler2008}.

In the context of statistical test of hypotheses, a \textit{compositional logic} is conveyed by an algebraic formalism that allows the evaluation of truth-functions of composite models and truth-values of composite hypotheses by algebraic operations on the corresponding truth functions of elementary models and truth values of elementary hypotheses. 
The e-value and the FBST have a rich, expressive, and intuitive compositional logic, while traditional truth-values and accompanying tests offered by either frequentist (classical) statistics, like the p-value, or by Bayesian statistics, like Bayes factors, have important and well-known deficiencies in this regard, specially in cases involving sharp statistical hypotheses. 
Furthermore, \textit{logically coherent} evaluations
and testing of sets and sub-sets of statistical hypotheses should render sequences of inferential reasoning that do not generate internal contradictions or anti-intuitive results. 
As expected, the e-value and the FBST comply with well-established rules of logical coherence, even in the case of sharp hypotheses, while traditional alternatives often fail to do so. 

Asymptotically, $\mbox{sev}(H|X)$ -- the standardized e-value -- shares several properties of the p-value, the well-known significance measure of frequentist statistics. 
This allows the use of the standardized e-value in frequentist-oriented applications, retaining, nevertheless, many theoretical characteristics of the Bayesian framework. 
Many theoretical developments and practical applications of the e-value and the FBST have already been published in the scientific literature, see \cite{PereiraStern2020}. 
Concise entries about the e-value and the FBST are available at the \textit{International Encyclopedia of Statistical Science} and online at Wiley's \textit{StatsRef}, see \cite{Diniz2022,Pereira2011}. 

Section 2 reviews the Bayesian statistical framework; 
Section 3 defines the e-value; 
Sections 4, 5, and 6 explain the invariance, asymptotic, and compositional properties of the e-value; 
Section 7 defines the GFBST - the Generalized Full Bayesian Significance Test and its logical properties; 
Sections 8 and 9 comment on computational implementation and give a detailed numerical example in model selection; 
Section 10 lists a representative assortment of articles from many practical applications of the e-value and the FBST already published in the scientific literature. 
Section 11 considers the philosophical consequences of the aforementioned developments by briefly commenting on the \textit{Objective cognitive constructivism} epistemological framework, which was specifically developed to accommodate the formal properties of the e-value and the FBST, and renders a naturalized approach to ontology and metaphysics. 
Section 12 presents some topics for further research at the interface between Logic and Statistics. 
One of the objectives of this paper is to foster work in areas of interface or overlap between Logic and Statistics and to stimulate greater cooperation between the two communities. Hence, for the sake of clarity and completeness, this paper includes some basic definitions and explanations in both areas that may be convenient to facilitate mutual understanding of more advanced topics. 

\section{Bayesian Framework} 

The basic definition of the e-value is given in the context of a standard parametric model in {Bayesian statistics}, where observable variables $x$ are generated with a probability density $p(x \g \theta)$. 
The observable (vectors of) random variables 
$x\in \mathscr{X} \subseteq \mathbb{R}^s$ belong to the model's \textit{sample space}, while the latent (non-observable) vector 
$\theta\in \Theta \subseteq \mathbb{R}^t$ belongs to the model's \textit{parameter space}. 
Statistical inference aims to acquire information about the unknown parameter $\theta$ from a sequence of observations, 
$X= [x^{(1)};x^{(2)};\ldots\, x^{(n)}]$, 
that for notational convenience are stacked in $n \times s$ matrix $X$. 

Since the true value of the parameter, $\theta^0$, is unknown, $\theta$ is treated in Bayesian statistics as a (vector) random variable. 
Available previous information about the parameter is represented by its \textit{a priori} density, $p_0(\theta)$. 
No previous information about $\theta$ is represented by a non-informative \textit{reference density}, $r(\theta)$. 
From \textit{Bayes rule}, it follows that, after having $n$ independent observations in dataset $X$, the available information about the parameter $\theta$ is represented by the \textit{posterior density} 
\[ p_n(\theta \g X) \ = \ (1/c_n) p_0(\theta) 
\pprod_{i=1}^n p(x^{(i)} \g \theta) \ = \ 
(1/c_n) p_0(\theta) p(X \g \theta) \ , 
\] 
where $c_n$ is the appropriate normalization constant. 

Whenever possible, we use a relaxed notation leaving implicit the condionalization on the observed dataset, for example, writing $p_n(\theta)$ instead of $p_n(\theta \g X)$. 
Further details about the Bayesian framework, including appropriate choices for non-informative reference densities, can be found in \cite{Box1973,DeGroot1970,DeGroot2011,Gelman2004,Jeffreys1961,Kadane2011,Kadane2016,Kapur1989,Stern2011a,Zellner1971}. 

A statistical hypothesis $H$ states that the parameter $\theta^0$ generating the observations $X$ belongs to the \textit{hypothesis' set}, $\Theta_H$, a region of 
$\Theta$ 
constrained by (vector) inequality and equality constraints, 
\[ \Theta_H = \{ \theta \in \Theta \g 
g(\theta) \leq \zero \wedge h(\theta) = \zero \} \ .
\] 
It is common practice in statistics to use a relaxed notation, writing $H$ for both the hypothesis statement and the hypothesis' set. 

The presence of $q$ equality constraint makes 
the hypothesis \textit{sharp} or \textit{precise}, namely, $H$ becomes a proper surface (sub-manifold) of dimension strictly lower than the dimension of the parameter space, that is, 
$h=\mbox{dim}(H) = t-q <\, t=\mbox{dim}(\Theta)$. 
In particular, a point hypothesis is a hypothesis of dimension $h=0$. 
Regarded as a subset of the $t$-dimensional parameter space, a sharp hypothesis has volume zero, and should therefore have zero posterior probability, at least for a regular (continuous and differentiable) posterior density, which is usually a natural assumption for statistical modeling, see \cite{Kelter2021b}. 
These conditions make $\Pr(H)$, the probability of hypothesis $H$, an inappropriate measure to 
evaluate sharp hypotheses, unless a special probability measure is created for the hypothesis set, a situation that may be computationally cumbersome and theoretically challenging, see \cite{SternPereira2014,Williams2021}. 
Nevertheless, sharp hypotheses have a prominent role in science, for the most important statements in exact sciences are \textit{natural laws}. These are formulated as equations in a theory of interest that, in turn, can be expressed by sharp hypotheses in statistical models used to empirically verify the same laws, see \cite{Esteves2019,Stern2011a,Stern2011b,Stern2017a,Stern2020}. 
Sharp hypotheses also naturally arise in legal applications, like auditing models for regulatory compliance and concerns related to the burden of proof in related legal cases, see \cite{Gaskins1992,Guy1994,Kokott1998,PereiraStern1999b,PereiraNS2000,Pigliucci2014,Stern2003,Stern2018a}. 
This state of affairs was the main motivation for defining the e-value as a Bayesian significance measure specially well-suited for sharp hypotheses. 

\section{The e-value} 
\label{eval}

$\ev(H\g X)\in [0,1]$, the e-value, or the epistemic value of hypothesis $H$ given the observed data $X$, or the evidence given by the observed data $X$ in favor of hypothesis $H$, and its complement, $\evb(H\g X)= 1-\ev(H\g X)$, are defined as follows: 
   
(i) $s(\theta)$, the \textit{surprise function} in a statistical model is defined as the quotient between the posterior and the reference densities in the statistical model, see \cite{Dubois1982,Evans1997,Good1983,Madruga2003,Royall1997,Shackle1968,Shackle1969,SternPereira2014}, 
\[ 
s(\theta)= {p_n(\theta)}/{r(\theta)} \ ; 
\] 

(ii) $s^*$, the maximum (or supremum) of the surprise function constrained to the hypothesis $H$, is defined as 
\[ s^* = \sup\nolimits_{\theta \in H} s(\theta) \ , 
\] 
A maximizing argument, 
$\theta^* \g s^* = s(\theta^*)$, 
is called a \textit{tangential point}, see Figure 1;

(iii) $T(v)$, the closed lower $v$-cut of the surprise function, see \cite{Dubois1982,SternPereira2014}, and its complement, the open upper $v$-cut of the surprise function, 
$\Tb(v)$, 
are defined as 
\[ 
T(v) = \{ \theta \in \Theta \g s(\theta) \leq v \}
\ , \ \ \ \ 
\Tb(v) = 
\{ \theta \in \Theta \g s(\theta) > v \} \ ; 
\] 
The upper $v$-cut at level $v=s^*$, $\Tb(s^*)$, is called the \textit{tangential set}, for its border corresponds to the contour line of the surprise function that is tangential to hypothesis $H$, see Figure 1.

(iv) $W(v)$, the \textit{truth function} or \textit{Wahrheitsfunktion} at level $v$, is defined as the posterior probability mass inside the lower $v$-cut of the surprise function, see \cite{Borges2007,SternPereira2014}, 
\[ 
W(v) = \int_{T(v)} p_n\left(\theta \right) d\theta \ , 
\] 
while its complement is defined as 
$\Wb(v) = 1-W(v)$; 

(v) $\ev(H\g X)$, the \textit{epistemic value} of hypothesis $H$ given the observed data $X$, 
is defined as the truth function $W(v)$ computed at level $v=s^*$, see \cite{Borges2007,PereiraStern1999,Pereira2008,SternPereira2014}, while its complement, $\evb(H\g X)$, 
the evidence given by the observed data $X$ against hypothesis $H$, has the complementary probability mass, 
\[ 
\ev(H\g X) = W(s^*) \ , \ \ \ \ 
\evb(H\g X) = \Wb(s^*) = 1-\ev(H)\ ; 
\] 
For further details and remarks on the development of the e-value and its historical predecessors, see \cite{PereiraStern1999,Stern2003,SternPereira2014}.

\setcounter{figure}{0} 
\begin{figure}[t] 
\centerline{ 
\includegraphics[height=3.0in, width=4.5in, 
trim= 0mm 0mm 0mm 0mm , clip]{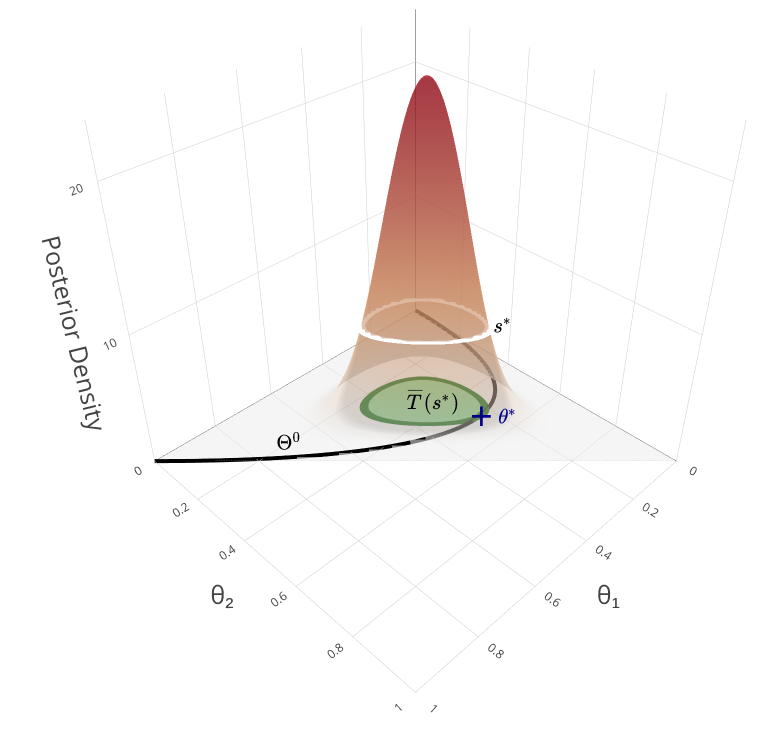} 
}

\caption{Flat prior, $r(\theta)\propto 1$; Surprise function equal to posterior density, $s(\theta)= p_n(\theta)$; Optimal $\theta^*\in H$; Contour curve at level $s^*=s(\theta^*)$; and Tangential highest probability density set, $\Tb(\theta^*)$.} 
\end{figure} 

In the special case of the flat reference density, $r(\theta) \propto 1$, the surprise function coincides with the posterior probability density function. Its upper cuts, $\Tb(v)$, coincide with \textit{highest probability density sets} (HPDS), making it easier to visualize all the elements defined in the last paragraph. 
Figure 1 depicts a simple example with a unimodal surprise function and a simply connected tangential set, see \cite[Sec.4.3]{PereiraStern1999}. However, in general, the posterior probability density function may have 
multiple local maxima, the tangential set may have multiple connected components, and all these components must be taken into account when computing the e-value. 
In this respect, the e-value differs from traditional uses of HPDS that only consider a single connected component, usually the connected component containing $\widehat{\theta}$, the unconstrained maximum a posteriory estimator, see \cite{Box1973,Dubois2004,Lindley1972,Salicone2007,SternPereira2014,West1997}. 

\section{Invariance and Reference Density} 
\label{reference} 
A \textit{reparameterization} of the parameter space, that is, a change in the coordinate system used to map $\Theta$, could ``stretch'' or ``compress'' the region surrounding a given location point in the map. 
This kind of effect is clearly visible comparing geographic maps using different cartographic projections. 
However, the probability mass inside a given region, as specified by a probability density function, should remain the same, regardless of the coordinate system in use. 
Hence, the value or ``height'' of a density function at a given point must change according to the coordinate system in use. 
This kind of reasoning explains the transformation rules specified by differential and integral calculus for modifying density functions according to transformations of the coordinate system, see \cite{Finetti1957,Flanders1974,Piskunov1969}. 
Nevertheless, since the surprise function is defined as the ratio of two densities, $s(\theta)=p_n(\theta)/r(\theta)$, its value remains unchanged by regular (continuous and differentiable) reparamaterizations of the parameter space, and the e-value inherits this important invariance property, see \cite{Borges2007,SternPereira2014}. 

\textit{Invariance} is a vital property for good statistical inference for, otherwise, conclusions reached by statistical analysis would depend on the coordinate system being employed, like units of measurement, map projections, forms of data visualization, and other idiosyncrasies or arbitrary choices. 
Invariance properties of the e-value can be further explained and analyzed in the context of \textit{Information geometry}, where the reference density, $r(\theta)$, can be seen as an implicit representation of the underlying \textit{Information metric} for the parameter space, see \cite{Amari1987,Amari2007,Amari2016,Basu1988,Bernardo2005,Borges2007,Eaton1989,Frieden2004,Gelman2004,Jeffreys1961,SternPereira2014}. 
\section{Asymptotic Consistency} 
\label{Asympt}

As the number $n$ of observations grows to infinity, that is, in the asymptotic limit $n\rightarrow \infty$, we expect a \textit{consistent} statistical procedure to reach the ``correct'' conclusion. 
This section analyses the asymptotic properties of the e-value that motivate its use to consistently evaluate a statistical hypothesis. 
It is easier to describe these asymptotic properties in terms of the \textit{standarized} version of the e-value, defined as follows: 

(i) $Q(d,z)$, the \textit{chi-square} cumulative distribution with $d\in \mathbb{N}_+$ degrees of freedom for random variable $z\in[0,\infty]$, is defined by the following analytical expression using the \textit{incomplete gamma function}: 
\[ 
Q(d,z) = 
\frac{\gamma(d/2, z/2)}{\gamma(d/2, \infty)} \ , 
\ \ \ \ 
\gamma(d,z) = \int_0^z y^{d-1}e^{-y} dy \ ; 
\] 

(ii) $\sigma(t,h,c)$, the \textit{standardization function} on arguments $t,h\in \mathbb{N}_+$ and $c\in [0,1]$, is defined by the expression 
\[ 
\sigma(t,h,c) = 
Q\left(t-h, Q^{-1}\left(t,c\right) \right) \ ; 
\] 

(iii) $\sev(H\g X)$, the standarized e-value of a hypothesis $H\subset \Theta$ of dimension 
$h=\dim(H) \leq t=\dim(\Theta)$, is defined as follows: 
\[ 
\sev(H\g X) =1-\sevb(H\g X) \ , 
\ \ \ \ 
\sevb(H\g X) = \sigma(t,h,\evb(H\g X)) \ ; 
\]

The standarized e-value has the following asymptotic properties, under usual continuity and differentiability regularity conditions, see \cite{Borges2007}: 

(a) If $H$ is true, i.e. $\theta^0\in H$; and 
$H$ is slack, i.e. $\dim(H)=\dim(\Theta)$; and 
$\theta^0$ is in the topological interior of $H$; 
Then $\ev(H\g X)$ converges (in probability) to 1, as the number of observations increases; 

(b) If $H$ is true, i.e. $\theta^0\in H$; and $H$ is sharp, i.e. $\dim(H)<\dim(\Theta)$; 
Then $\sev(H\g X)$ converges (in distribution) to $U[0,1]$, the uniform distribution in the unit interval. 
The standardized e-value shares these asymptotic properties with the p-value -- the well-known and widely used significance measure of frequentist statistics. 
Hence, in a frequentist-oriented application, one could use the e-value as a Bayesian replacement for the p-value, preserving already familiar decision procedures and their interpretations. 
Moreover, in several modeling applications, the e-value exhibits better convergence characteristics than those of the p-value, see \cite{Bernardo2012,Lauretto2003,Lauretto2007,Lauretto2005a,Lauretto2005b}. 
Further consistency properties of the e-value and higher order asymptotic approximations have been studied and developed in \cite{Cabras2015,Pinto2012,Ranzato2018,Ruli2020,Ruli2020b,Ventura2013,Ventura2014,Ventura2016}.

\section{Compositional Logic of e-values} 
The e-value and its truth function yield simple algebraic expressions that facilitate the study of a composite model, $M$, build by serial coupling of $k$ independent statistical models,   
indexed on $j=1\ldots k$, 
$M^{(j)}=\{\Theta^{(j)}, p_0^{(j)}, p_n^{(j)}, r^{(j)} \}$.  
In this setting, $q$ alternative (or parallel) sets of hypotheses, $H^{(i,j)}$, indexed on $i=1\ldots q$, are provided for evaluation. 
This situation is somewhat analogous to the analysis of parallel-serial systems in reliability theory, see \cite{Barlow1981,Borges1979,Kaufmann1977} and Figure 2. 
Let us first consider a pure serial system, consisting of $k$ individual models, $M^{(j)}, j=1\ldots k$, each one contemplating a single hypothesis, $H^{(1,j)}$. 

Since the individual models are independent, the joint posterior density, reference, and surprise function of the composite model are the product of the corresponding individual functions, that is, 
\[
p_n(\theta) \ = \ \pprod_{j=1}^k p_n^{(j)}(\theta^{(j)}) 
\ , \ \ \ 
r(\theta) \ = \ \pprod_{j=1}^k r^{(j)}(\theta^{(j)}) \ , 
\] 
\[ 
s(\theta) \ = \ \pprod_{j=1}^k s^{(j)}(\theta^{(j)}) 
\ , \ \ 
\theta = [\theta^{(j)}, \ldots\, \theta^{(k)}] \in \Theta \ .
\]

The Mellin convolution of the cumulative distributions for (scalar) random variables $x$ and $y$, $F(x)\otimes G(y)$, gives the cumulative distribution of the product $z=xy$. 
Hence, the composite model's truth function is the \textit{Mellin convolution} of the individual truth functions, that is,
\[ W(v) = 
\bigotimes_{1\le j\le k}W^{(j)} \ = \ 
W^{(1)}\otimes W^{(2)}\ldots\, \otimes W^{(k)}(v) \ . 
\] 
Further properties and interpretations of the Mellin convolution, and detailed analytical procedures and numerical algorithms for its efficient computation, can be found in \cite{Bracewell1986,Kaplan1987,Springer1979,Williamson1989}. 
Matlab code for computing the composite model's truth function in the discrete case, i.e. for step ladder cumulative functions, is presented in \cite{Borges2006}.   
 
\setcounter{figure}{1} 
\begin{figure}[t] 
\centerline{ 
\includegraphics[height=1.5in, width=3.5in, 
trim= 0mm 0mm 0mm 0mm , clip]{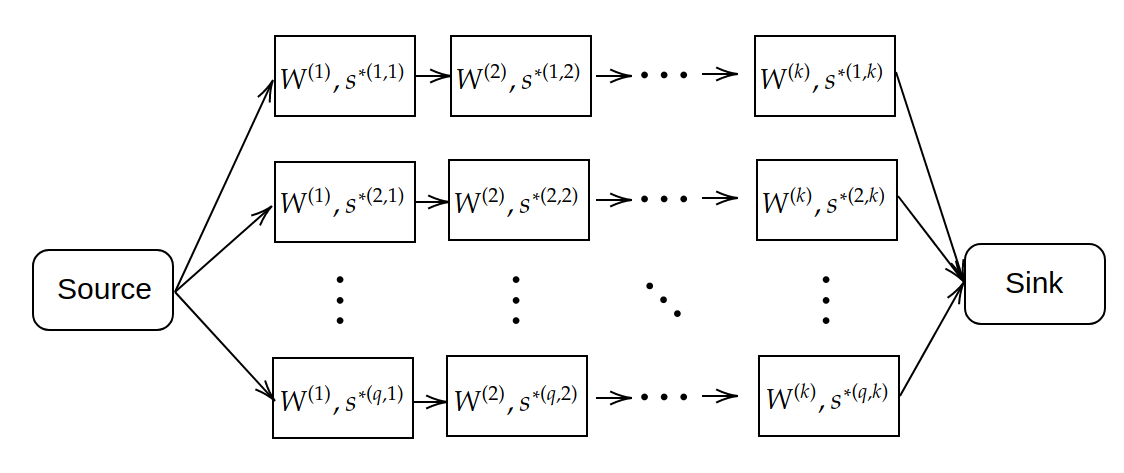} 
}

\caption{Diagrammatic representation of parallel hypotheses in serial models composition of a composite hypothesis, $H$, whose e-value, $\ev(H)$, can be computed from the truth-functions, $W^{(j)}$, of each single model, $M^{(j)}$, and the optimal surprise values, $s^{*(i,j)}$, of each single hypothesis, $H^{(i,j)}$.} 
\end{figure} 

Next, let us consider a composite hypothesis $H$ expressed in \textit{conjunctive form}, where the conjunction symbol ($\wedge$) stands for the \textit{and} logical operator, 
\[ 
H = \bbigwedge_{j=1}^k H^{(1,j)} = 
H^{(1,1)} \wedge H^{(1,2)} \ldots\, \wedge H^{(1,k)} \ . 
\] 

In the conjunctive hypothesis, the surprise 
function attains a maximum (or supremum) value equal to the product of individual maxima (or suprema). Therefore, the e-value of the conjunctive composite hypothesis is given by 
\[
\ev(H) \ = \ W(s^*) \ = \ 
\bigotimes_{1\le j\le k} 
W^{(j)} \left( \pprod_{j=1}^k s^{*(1,j)} \right) \ .
\]

Finally, let us consider a composite hypothesis $H$ expressed in \textit{disjunctive normal form}, where the disjunction symbol ($\vee$) stands for the \textit{or} logical operator, 
\[
H \ = \ \bbigvee_{i=1}^q \bbigwedge_{j=1}^{k} H^{(i,j)} = 
\left( \bbigwedge_{j=1}^{k} H^{(1,j)} \right) 
\vee 
\ldots\, \vee 
\left( \bbigwedge_{j=1}^{k} H^{(q,j)} \right). 
\]
The corresponding e-value requires the conjunctive surprise to be maximized over the finite set of $q$ disjunctive alternatives, 
\[ 
\ev(H) \ = \ 
\ev \left( \bbigvee_{i=1}^{q} 
\bbigwedge_{j=1}^{k} H^{(i,j)} \right) 
\ = \ 
\] 
\[
W(s^*) \ = \ 
\mmax_{i=1}^{q} 
\bigotimes_{1\le j\le k} W^{(j)} 
\left( \pprod_{j=1}^{k} s^{*(i,j)} \right) \ . 
\]

Interestingly, in the last expression, if all elementary hypotheses have either null or full epistemic value, that is, if, 
$\forall (i,j), \ev(H^{(i,j)})\in \{0,1\}$, the evaluation of $\ev(H)$ reduces to the the evaluation in classical logic of a corresponding expression preserving the same normal structure \cite{Borges2007}. 
Another interesting aspect of the last expression is how conjunction operators are translated into product operations, while disjunction operators are translated into maximization operations. 
This kind of translation characterizes the e-value as a \textit{possibilistic} \textit{abstract belief calculus}. 
For further details and comments, see \cite{Borges2007,Darwiche1992,Darwiche1993,Dubois1982,Dubois2012,Klir1988,Stern2003,Stern2004,SternPereira2014,Thulin2014,Vieland2019}.

\section{Statistical Hypotheses Tests and the GFBST}

In the practice of statistics, it is often necessary to make a decision to either \textit{reject}, or to remain \textit{undecided}, or to \textit{accept} a hypothesis $H$ (in statistics, the alternative undecided is also called \textit{agnostic}). 
In \textit{multi-valued logics}, these three alternatives are traditionally encoded by the logical values $\{0, \text{1/2}, 1\}$, see \cite{Epstein1993,Epstein1974,Gottenwald2001,Rine1984}. 
In \textit{modal logics}, these three alternative attitudes consider the hypothesis to be \textit{impossible}, \textit{contingent}, or \textit{necessary}, see \cite{Carnielli2008a,Garson2006,Priest2001}. 
Table 1 presents the corresponding \textit{modal operators} together with some relevant compositions obtained by negation ($\lnot$, the \textit{not} operator), disjunction ($\vee$, the \textit{or} operator) and conjunction ($\wedge$, the \textit{and} operator).

\setcounter{figure}{2} 
\begin{figure}[t] 
\centerline{ 
\includegraphics[height=2.5in, width=2.5in, 
trim= 0mm 0mm 0mm 0mm , clip]{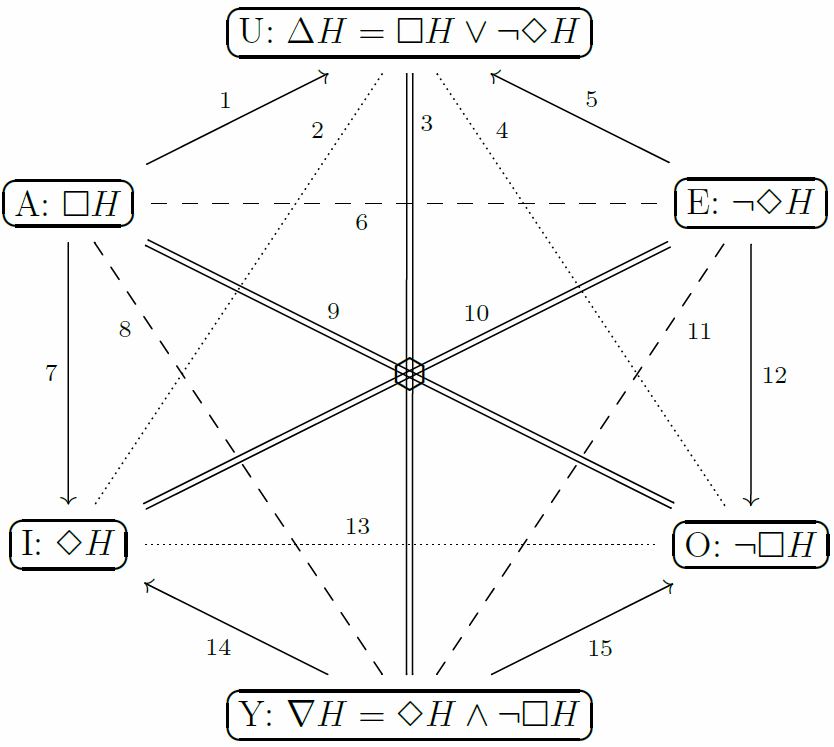} 
}

\mbox{} \vspace{0mm} \mbox{} \\ 


\begin{centering} 
\noindent 
{Figure 3: Hexagon of opposition depicting modal operators of necessity, $\square$; possibility, $\Diamond$; contingency, $\nabla$; and non-contingency, $\Delta$; some compositions by negation, $\lnot$; disjunction, $\vee$; and conjunction, $\wedge$; and logical relations of implication, $\rightarrow$; contrariety $-\,-$; sub-contrariety, $\cdots$; and contradiction, $=\!\!=$.} 
\end{centering} 

\mbox{} \vspace{0mm} \mbox{} \\ 

\centering 
\begin{tabular}{|l|c|c|l|}
\hline
Name & Modality & Equivalence & Interpretation \\
\hline
Necessity & $\square H$ & $\Delta H \wedge \Diamond H$ 
& $H$ is accepted \\
Impossibility & $\lnot \Diamond H$ & 
$\Delta H \wedge \lnot \square H$ & $H$ is rejected \\
Contingency & $\nabla H$ & 
$\Diamond H \wedge \lnot \square H$ & $H$ is not decided \\
Possibility & $\Diamond H$ & $\square H \vee \nabla H$ & $H$ is not rejected \\
Non-necessity & $\lnot \square H$ & $\lnot \Diamond H \vee \nabla H$ & $H$ is not accepted \\
Non-contingency & $\Delta H$ & $\square H \vee \lnot \Diamond H$ & $H$ is decided \\
\hline
\end{tabular}

\mbox{} \vspace{2mm} \mbox{} \\ 

\centerline{Table 1:Basic modalities of agnostic hypothesis tests.} 
\end{figure}

The modalities in Table 1 are related by logical relations depicted in diagrammatic form by the \textit{hexagon of opposition} in Figure 3, where: 
Arrows ($\longrightarrow$) indicate logical \textit{implication}; 
Dashed-lines ($---$) indicate \textit{contrariety}, connecting modalities that cannot both be true; 
Dotted-lines ($\cdots\cdots$) indicate \textit{sub-contrariety}, connecting modalities that cannot both be false; and 
Double-lines ($=\!\!=$): indicate \textit{contradiction}, connecting modalities that can neither both be true nor both be false. 
The vowels at the vertices of the hexagon are historical labels inherited from medieval logic, see \cite{Dumitriu1977}.

The logical relations depicted in the hexagon of opposition are considered basic principles of rational argumentation, intuitive to, and tacitly assumed by most users and, sometimes, even soft coded in popular culture or hard coded in natural language, see \cite{Beziau2012,Blanche1966,Dumitriu1977,Gallais1982,Izbicki2015,Stern2022,Stern2023}. 
Therefore, departing from these principles may make arguments counterintuitive, create barriers to the natural flow of communication, or even lead to misunderstandings. 
This situation justifies the development of inference methods and supports ways of reasoning that are fully compliant with these basic principles.

A \textit{statistical test of hypotheses} is a statistical procedure used to make a required decision to either reject, or remain undecided, or accept a hypothesis $H$. 
It is convenient to use a statistical test based on an already available significance measure, $\mu(H)\in [0,1]$. 
Such a statistical test can be regarded as a discretization map, $\delta(\mu(H))$, that collapses the interval $[0,1]$ into the ternary set $\{0, 1/2,1\}$. 
Some applications may only allow decisions in a binary subset, like $\{0, 1/2\}$, although a full logical analysis of these tests requires (at least) a ternary decision space, see \cite{Esteves2023,Izbicki2015,Izbicki2022,Silva2015,Stern2018c}. 
  
The remainder of this section examines logical consistency conditions for such tests, and provides further justification for the use of the e-value as an appropriate significance measure. 
The preceding sections paid special attention to sharp hypotheses. In contrast, this section looks at sharp or slack hypotheses with equal interest. 
Although motivated by the need to evaluate sharp hypotheses, we should remark that nothing prevents the e-value from being computed for slack hypotheses. 

As already examined in Section 2, by definition, a statistical hypothesis $H$ states that the true parameter value belongs to the hypothesis set, that is, that $\theta^0\in H$. 
From this definition, in conjunction with classical principles or rational argumentation encoded in the hexagon of opposition, we now present three sets of conditions for logical consistency that ought to be required from a rational hypothesis test, namely, \textit{invertibility}, \textit{monotonicity}, and \textit{consonance}, see \cite{Izbicki2015,Esteves2016,Esteves2019,Stern2018c}. 

\mbox{} \vspace{-2mm} \mbox{} 
   
\textit{Invertibility:} 
Every hypothesis $H\subset \Theta$ automatically defines its complement, $\Hb = \Theta -H$. 
Hence, the true value of the parameter, 
$\theta^0\in \Theta$, must either belong to $H$ or belong to its complement. 
Therefore, the following implications ought to hold: 
\begin{description} 
\setlength\itemsep{0mm} 

\item[(I.i)] \textit{Necessity inversion:} 
$H$ is necessary if and only if $\Hb$ is impossible, that is, 
$\square H \Leftrightarrow \lnot \Diamond \Hb$; 

\item[(I.ii)] \textit{Possibility inversion:} 
$H$ is possible if and only if $\Hb$ is unnecessary, that is, 
$\Diamond H \Leftrightarrow \lnot \square \Hb$;
and 

\item[(I.iii)] \textit{Contingency inversion:} 
$H$ is contingent if and only if $\Hb$ is contingent, that is, 
$\nabla H \Leftrightarrow \nabla \Hb$. 
\end{description} 
Figure 4 (left) gives a diagrammatic representation of invertibility relations (and their complements obtained by negation) using fragments of the hexagon of opposition. 
  
\mbox{} \vspace{-2mm} \mbox{} 
  
\textit{Monotonicity:} 
Let us consider a larger hypothesis $H'$ that contains the original hypothesis $H$, that is $H'\supset H$. 
In this setting, the following implications ought to hold: 
\begin{description} 
\setlength\itemsep{0mm} 
\item[(M.i)] \textit{Monotonic necessity:} If $H$ is necessary, so must be $H'$, that is, $\square H \Rightarrow \square H'$; and 

\item[(M.ii)] \textit{Monotonic possibility:} If $H$ is possible, so must be $H'$, that is, 
$\Diamond H \Rightarrow \Diamond H'$. 
\end{description} 
Figure 4 (right) gives a diagrammatic representation of these monotonicity relations (and their complements obtained by negation) using fragments of the hexagon of opposition. 
  
\mbox{} \vspace{-2mm} \mbox{} 
  
\textit{Consonance:}
Let $H_i$, for $i\in I$, be a collection of hypotheses referenced by the index set $I$. 
The following implications ought to hold: 
\begin{description} 
\setlength\itemsep{0mm} 
\item[(C.i)] \textit{Union consonance}: If the \textit{join}  hypothesis, made by the union of all indexed hypotheses, is possible, then at least one of the joining hypotheses is possible, that is, 
$$ \Diamond ( \cup_{i \in I} H_i) \Rightarrow 
\exists i \in I \g \Diamond H_i \ . $$ 
\item[(C.ii)] \textit{Intersection consonance}: 
If the \textit{meet} hypothesis, made by the intersection of all indexed hypotheses, is unnecessary, then at least one of the meeting hypotheses is unnecessary, that is, 
$$ \lnot \square (\cap_{i \in I} H_i) \Rightarrow 
\exists i \in I \g \lnot \square H_i 
\ \ \Leftrightarrow \ \ 
\forall i \in I , \ \square H_i \Rightarrow 
\square (\cap_{i \in I} H_i) \ .$$ 
\end{description} 
Figure 5 gives a diagrammatic representation of consonance relations for a collection of three hypotheses, $\{A, B, C\}$.

\setcounter{figure}{3} 
\begin{figure}[t] 
\centerline{ 
\includegraphics[height=1.5in, width=2.0in, 
trim= 0mm 0mm 0mm 0mm , clip]{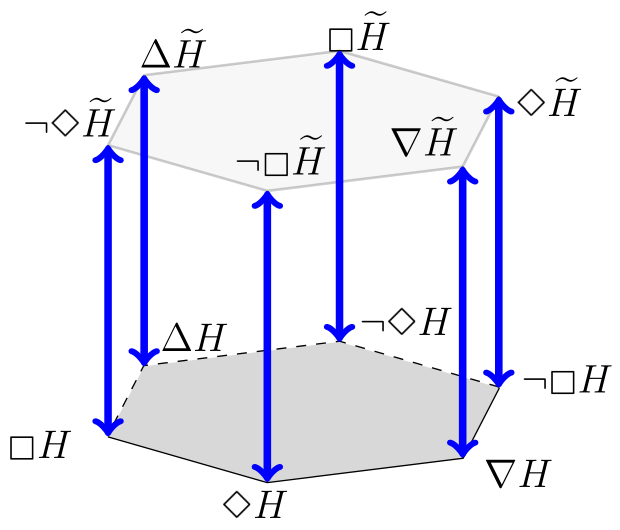} 
\mbox{} \hspace{2mm} \mbox{} 
\includegraphics[height=1.5in, width=2.0in, 
trim= 0mm 0mm 0mm 0mm , clip]{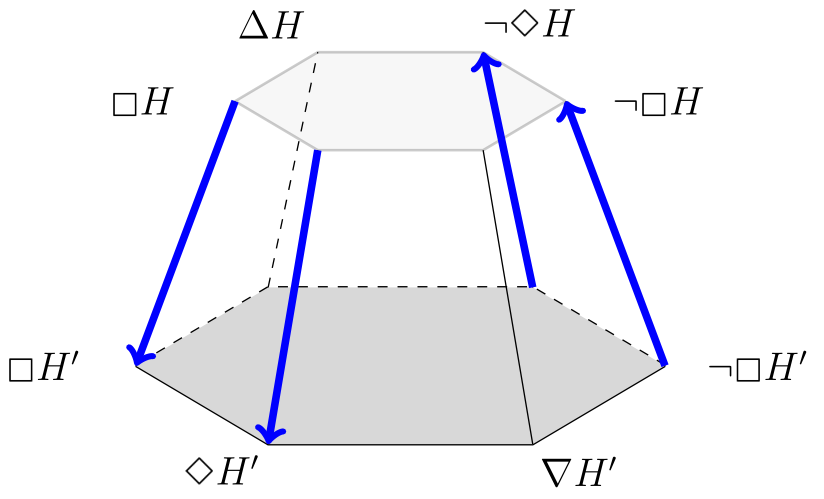} 
}

\caption{Invertibility (left) and Monotonicity (right) logical relations between statistical modalities applied to a hypothesis, $H$, and its complement, $\Hb$. } 
\end{figure}

\setcounter{figure}{4} 
\begin{figure}[t] 
\centerline{ 
\includegraphics[height=1.5in, width=4.0in, 
trim= 0mm 0mm 0mm 0mm , clip]{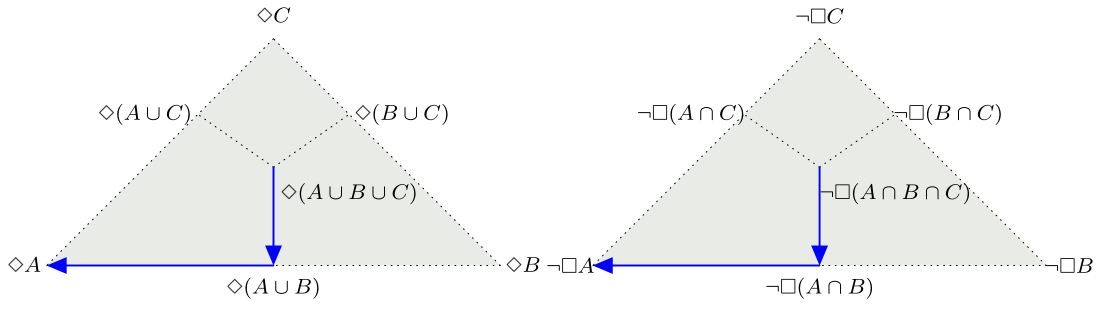} 
}

\caption{Consonance logical relations for the union or intersection of three hypotheses.} 
\end{figure}

\setcounter{figure}{5} 
\begin{figure}[h] 
\centerline{ 
\includegraphics[height=1.3in, width=4.0in, 
trim= 0mm 0mm 0mm 0mm , clip]{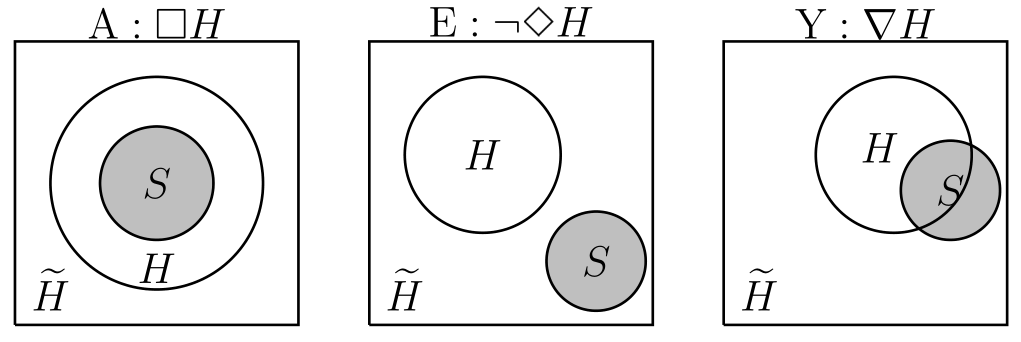} 
}

\caption{Testing $H$ with region estimator $S$: Accept if $S\subset H$, reject if $S\subset \Hb$, undecided otherwise. 
} 
\end{figure}

Arguments of analysis lead to a characterization of test procedures that comply with the aforementioned logical conditions of monotonicity, invertibility, and consonance, see \cite{Esteves2016,Izbicki2022}. 
These tests must be based on a \textit{region estimator}, $S\subset \Theta$, for the location of the true parameter value, $\theta^0$. 
Such a \textit{region estimator test} accepts $H$ if $S\subset H$, rejects $H$ if $S\subset \Hb$, and remains undecided if $S$ properly intersects both $H$ and $ \Hb$, as depicted in Figure 6.

As mentioned in the beginning of this section, it is convenient to use a statistical test, 
$\delta(\mu(H))\in \{0,1/2,1\}$, that is based on an already available significance measure, 
$\mu(H)\in [0,1]$. 
Nevertheless, the following logical compatibility condition, between the statistical test and its underlying significance measure, ought to hold: 
\[ 
\mu(H^{(1)})\geq \mu(H^{(2)}) \Rightarrow
\delta(\mu(H^{(1)}))\geq \delta(\mu(H^{(2)})) \ . 
\] 
This logical compatibility condition is inspired by several legal principles, like the ones known in the juridical literature as \textit{onus probandi}, 
\textit{in dubio pro reo}, and \textit{principle of proportionality}, see \cite{Barak2012,Engle2012,Gaskins1992,Guy1994,Hulsroj2013,Kokott1998,Pigliucci2014,Stern2003,Romero1991,Zeleny1982}. 
The \textit{onus probandi} and \textit{in dubio pro reo} principles can be interpreted as requiring the monotonicity and consonance properties examined in this section, while the \textit{principle of proportionality} 
can be interpreted as follows: 
-- Judging the acceptability of hypotheses $H^{(1)}$ and $H^{(2)}$, if the former has a better grade (a higher significance measure) than the latter, then fair judgments must render verdicts (decisions) that are as least as favorable to the former than to the latter. 

From the definition of the e-value, it is clear that an upper $v$-cut $\Tb(v)$ can be used as a region estimator for $\theta^0$. 
Accordingly, the corresponding region estimator test rejects $H$ if its e-value stays below the threshold $c=W(v)$, that is, $H$ is rejected if $\ev(H)<c$. 
By invertibility, $H$ is accepted if $\ev(\Hb)<c$. 
This is the GFBST - the \textit{Generalized Full Bayesian Significance Test}, see 
\cite{Esteves2016,Fossaluza2017,Izbicki2015,Madruga2001,Silva2015,Silva2018,Stern2018c}. 
One can check that the GFBST generates logical modalities that are compatible with all logical requirements examined in this section, and also that it is logically compatible with its underlying significance measure, the e-value. 
Moreover, one can also check that all good invariance and asymptotic properties of the e-value are inherited by the GFBST. 
Furthermore, it is interesting to remark that, under usual regularity conditions (continuity and differentiability), if $\ev(H)<1$, a sharp hypothesis $H$ may either remain undecided or be rejected by the GFBST, but never be accepted. 
Finally, for simplicity and limitations of space, the exposition given in this article first defines the e-value and the GFBST, and then derives or explains their logical properties. 
In contrast, in \cite{Esteves2016,Esteves2019,Izbicki2015,Izbicki2022,Silva2015,Stern2018c}, the authors show that it is possible to go the other way around, demonstrating that the aforementioned logical properties render a complete characterization of the GFBST as a region estimator test.

\section{Computational Implementation} 

Numerical computation of the e-value can be accomplished in two basic steps, namely, 
\begin{description} 

\item[\textit{Integration step}] -- computing an approximation of the truth function $W(v)$. 
In most cases, this task can be accomplished by using a computational condensation procedure, like \cite{Kaplan1987,Williamson1989}, on a numerical sequence, $s^{(k)}=s(\theta^{(k)})$, obtained by sampling from the posterior distribution, $p_n(\theta)$.  
Sampling sequences can, in turn, be generated using Markov Chain Monte-Carlo methods, see   \cite{Brooks2011,Gilks1996,Hammersley1964}, or variations thereof, like the Hit-and-Run algorithm, see  \cite{Chen2000,Karawatzki2005}; 

\item[\textit{Optimization step}] -- finding the optimum of the objective function $s(\theta)$ under the constraints imposed by the hypothesis. 
In most cases, this task can be accomplished by standard methods of continuous constrained optimization, like Generalized Reduced Gradient, Sequential Quadratic Programming, Generalized Augmented Lagrangian, or Proximal Point methods, see 
\cite{Birgin2014,Fang1997,Kapur1992,Luenberger2008,Minoux1986}. 
The following numerical optimization solvers 
are readily available and potentially applicable: \\ 
\textit{TANGO} -- Trustable Algorithms for Nonlinear General Optimization, a joint project from USP -- the Universtity of S\~{a}o Paulo, and UNICAMP -- the University of Campinas, that offers excellent code under General Public License (GNU); and \\ 
\textit{Gurobi} -- a robust, high-performance, and versatile optimization software that can currently be licensed for free academic use. \\ 
Some optimization problems with many local maxima can be solved indirectly by repeated local optimization from candidate starting points selected at the integration step, or directly by stochastic optimization methods based on Simulated Annealing, see \cite{Pflug1996,Spall2003,Stern1992,Voigt1990,Wah2008}. 

\end{description}

\section{Numerical Example in Model Selection}

Figure 7 (left) depicts the polynomial fitting problem for the classical Sakamoto et al. benchmark dataset presented in \cite[ch.8]{Sakamoto1986}. 
This dataset, given in Table 1, was generated by a simulation from the i.i.d. (independent and identically distributed) stochastic process 
\[ 
y_i = g(x_i) +N(0\,,\, 0.1^2) \ , \ \ 
g(x) = \exp((x -0.3)^2) -1 \ , 
\] 
that is, the points $y_i$ were generated at the grid $x_i= (i-1)0.05$, for $i=1,\ldots 21$, by adding to the exponential \textit{target function}, $g(x)$, a Gaussian random noise with mean $\mu=0$ and standard deviation of $\sigma=0.1$.

The \textit{linear regression} polynomial model of order $k$ explains vector $y$ in the dataset, 
using as explanatory variables integer powers up to order $k$ of the grid points in vector $x$, plus an i.i.d. Gaussian random noise with standard deviation $\sigma$, that is, 
\[ 
y= \beta_0 x^0 +\beta_1x^1 +\beta_2x^2 \ldots\, +\beta_kx^k +N(0,\sigma I) \ . 
\] 
Using the weakly informative prior density, 
$p_0(\beta,\sigma)= 1/\sigma$, 
the posterior density for this statistical model can be conveniently written, see \cite{DeGroot1970,Hocking1985,Zellner1971}, 
in terms of the matrix of explanatory variables, 
$X=[x^0;x^1;\ldots\, x^k]$, 
the maximum a posteriori estimators of the model's parameters, 
$\hat{\beta} = (X^tX)^{-1}X^ty$\,, 
and the auxiliary quantities \ 
$\hat{y} = X\hat{\beta}$ \ and \ 
$s^2 = (y -\hat{y})^t(y -\hat{y})/(n-k) \ ,$ 
\[ 
p_n(\beta,\sigma\g y,x^0,\ldots\, x^k) = 
\frac{1}{\sigma^{n+1}} \exp 
\left( 
\frac{-1}{2\sigma^2} \left( (n-k)s^2 +(\beta -\hat{\beta})^tX^tX(\beta -\hat{\beta}) 
\right) \right) \, . 
\]

The model selection problem that naturally arises in this context can be stated as the question: 
From a statistical point of view, which order $k$ gives the best polynomial fit for this dataset? 

The first column of Table 3 presents the quadratic norm \textit{empirical error}, 
$R_{emp} = || \hat{y} -y ||_2$\,, 
of each model, see \cite{DeGroot1970,Hocking1985,Zellner1971} for further details. 
Polynomials of increasing order have increasing flexibility, making $R_{emp}$ monotonically decrease until it reaches zero for order $k=n-1$, when the model becomes just an interpolating polynomial for the dataset.

\setcounter{figure}{6} 
\begin{figure}[h] 
\centerline{ 
\includegraphics[height=2.5in, width=4.4in, 
trim= 0mm 0mm 0mm 0mm , clip]{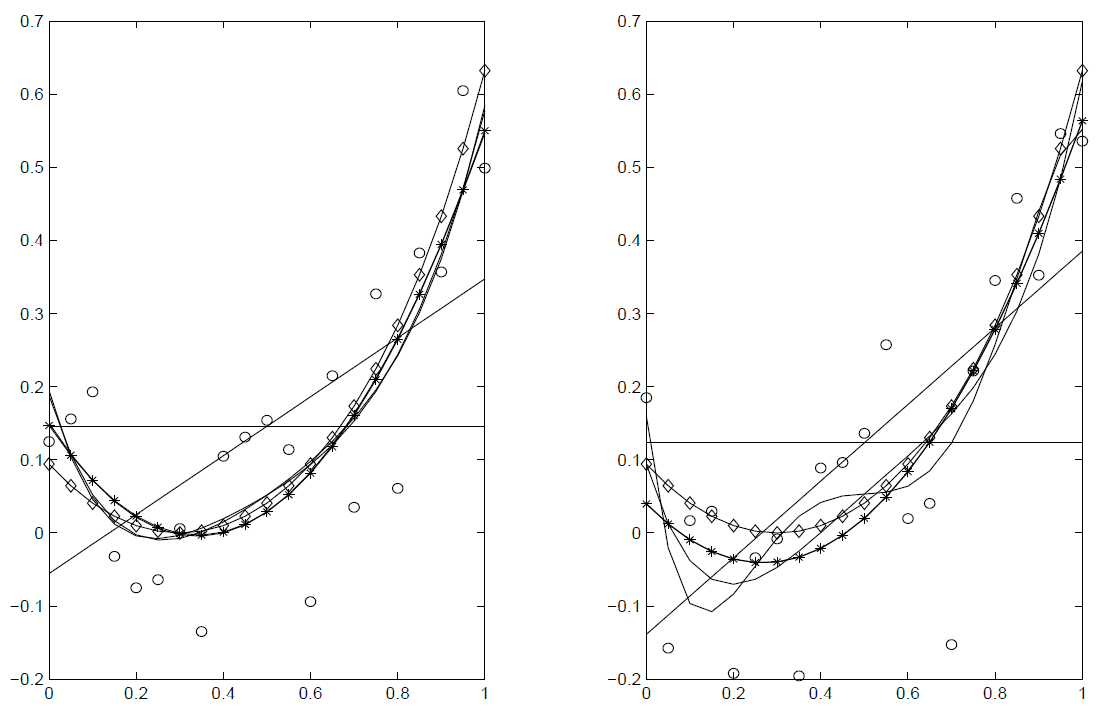} 
}
\caption{Benchmark (left) and alternative (right) data points, $\circ$, for Sakamoto et al. \cite[ch.8]{Sakamoto1986} polynomial fitting problem; 
Exponential target function, $\diamond$; 
Best fitted polynomials, --- , of order 0 to 4;
and Order 2 polynomial, $*$, that renders the smallest regularized error for both datasets.} 
\end{figure} 


\begin{table}[h!] 
\begin{center} 
\mbox{} \vspace{2mm} \mbox{} \\ 
\begin{small} 
\begin{tabular}{|r r r|r r r|r r r|} 
\hline 
$i$ & $x_i$ & $y_i$ & $i$ & $x_i$ & $y_i$ & $i$ & $x_i$ & $y_i$ \\ 
\hline 
1 & 0.00 & 0.125 & 8 & 0.35 & -0.135 & 15 & 0.70 & 0.035 \\ 
2 & 0.05 & 0.156 & 9 & 0.40 & 0.105 & 16 & 0.75 & 0.327 \\ 
3 & 0.10 & 0.193 & 10 & 0.45 & 0.131 & 17 & 0.80 & 0.061 \\ 
4 & 0.15 & -0.032 & 11 & 0.50 & 0.154 & 18 & 0.85 & 0.383 \\ 
5 & 0.20 & -0.075 & 12 & 0.55 & 0.114 & 19 & 0.90 & 0.357 \\ 
6 & 0.25 & -0.064 & 13 & 0.60 & -0.094 & 20 & 0.95 & 0.605 \\ 
7 & 0.30 & 0.006 & 14 & 0.65 & 0.215 & 21 & 1.00 & 0.499 \\ 
\hline 
\end{tabular}
\mbox{} \vspace{4mm} \mbox{} \\ 
\begin{tabular}{|c|c|c|c|c|c|c|c|c|} 
\hline 
Order & $R_{EMP}$ & $R_{FPE}$ & $R_{SBC}$ & $R_{GCV}$ & $R_{SMS}$ & $R_{AIC}$ & $\ev(H)$ & $\sev(H)$ \\ 
\hline 
0 & 0.03712 & 0.04494 & 0.04307 & 0.04535 & 0.04419 & -07.25 & 0.000 & 0.000 \\ 
1 & 0.02223 & 0.02964 & 0.02787 & 0.03025 & 0.02858 & -20.35 & 0.009 & 0.001 \\ 
2 & 0.01130 & 0.01661 & 0.01534 & 0.01724 & 0.01560 & -32.13 & 0.013 & 0.000 \\ 
3 & 0.01129 & 0.01835 & 0.01667 & 0.01946 & 0.01667 & -30.80 & 0.999 & 0.818 \\ 
4 & 0.01088 & 0.01959 & 0.01751 & 0.02133 & 0.01710 & -29.79 & 0.995 & 0.403 \\ 
5 & 0.01087 & 0.02173 & 0.01913 & 0.02445 & 0.01811 & -27.86 & 0.999 & 0.641 \\ 
\hline 
\end{tabular} 
\end{small}
\mbox{} \vspace{4mm} \mbox{} \\ 
\centerline{Table 2: Sakamoto et al. \cite[ch.8]{Sakamoto1986} benchmark dataset for the polynomial fitting problem} \vspace{1mm} 
\centerline{Table 3: Model selection by regularized empirical errors and $\ev(\beta_k=0 \g y,\, [x^0;x^1;\ldots\, x^k])$ } 
\end{center} 
\end{table} 

Nevertheless, the predictive power of these models 
does not monotonically increase with the decrease of empirical error. 
From Figure 7 (right), it is easy to see that, if the order of the polynomial is too high, the model becomes excessively complex and over-fits the data, resulting in poor generalizations (predictions) outside the grid points. 

Columns 2 to 6 in Table 3 present a variety of penalized errors, 
$R_{pen} = r(d,n) R_{emp}$\,, 
where the regularization factor, $r(d,n)$, 
penalizes the statistical model complexity. 
Specifically, the factor, $r(d,n)$, increases with the dimension of the model's parametric space, $d$, relative to the number of available data points, $n$. 
For the linear regressions at hand, having parameters $\sigma$ and $\beta_0, \beta_1,\ldots\, \beta_k$, $d=k+2$. 
The penalized or regularized errors in Table 3 are defined by the following regularization factors, based on the quotient $q=(d/n)$. 
The corresponding \textit{model selection criteria} chose the model with smallest regularized error; see \cite{Cherkassky1998,Sakamoto1986} for theoretical or heuristic justifications and further details.  
\begin{itemize} 
\setlength\itemsep{-1mm} 

\item Akaike's final prediction error, 
$r_{FPE}= (1+q)/(1-q) \, ;$ 

\item Schartz' Bayesian criterion, 
$r_{SBC}= 1 +\ln(n)q/(2-2q) \, ;$ 

\item Generalized cross validation, 
$r_{GCV}=(1-q)^{-2} \, ;$ 

\item Shibata model selector 
$r_{SMS}=1+2q \, .$ 
\end{itemize} 

The last columns in Table 3 present, in the context of a linear regression polynomial model of order $k$, 
the e-value and its standarized version for the statistical hypotheses stating that the parameter of higher order is null, that is, 
$\ev(H)=\ev(\beta_k=0 \g y, [x^0;x^1;\ldots\, x^n])$ and $\sev(H)$.  
This example demonstrates that the e-value and the FBST can replace traditional model selection criteria based on penalized prediction errors. 
Moreover, using a selection criterion based on the FBST renders a decision process that is invariant, asymptotically consistent, and logically coherent, as explained in the previous sections. 
In contrast, selection criteria based on penalized prediction errors lack even the most basic invariance properties, even though alternative model selection criteria exist that are invariant and have other desirable statistical properties, see for example \cite{Hutter2005,Hutter2008}. 
Furthermore, the use of the e-value and the FBST for model selection can be extended to a variety of non-nested (separate) or nested families of hypotheses, including Bayesian classifiers, as analyzed in \cite{Bernardo2012,Camargo2012,Hubert2009,Lauretto2003,Lauretto2005a,Lauretto2005b,Lauretto2007,Maranhao2012}.

\section{Applications} 

The development of statistical significance measures and tests may be motivated by their intended theoretical properties, which, in turn, may be inspired by epistemological desiderata. 
Nevertheless, these significance measures and tests must also prove themselves on the battlefields of science and technology as effective, efficient, robust, and reliable tools for the trade. 
A collection of over a hundred published applications of the e-value has been compiled in the 2020 survey \cite{PereiraStern2020}. 
This section gives a selection of references from this survey, organized by application area.

\begin{itemize} 
\setlength\itemsep{-1mm} 

\item 
Testing covariance structures in multivariate Normal models, treating in a unified way several alternative hypotheses (often treated as special cases in the literature): \cite{Lauretto2003,Vikas2016}; 
\item 
Testing unit root and cointegration hypotheses in time series, using plain and simple forms of prior information like flat or Jeffreys priors (no need for artificial priors): \cite{Chen2015,Diniz2011,Diniz2012a,Vosseler2016}; 
\item 
Solving Bayesian classification problems and testing nested and non-nested or separate hypotheses: \cite{Assane2018,Assane2019,Lauretto2005a,Lauretto2005b,Lauretto2007,Pereira2b}; 
\item 
Analyzing systems' reliability from failure datasets: 
\cite{Irony2002,Maranhao2012,Rodrigues2006} 

\item 
Testing dependence structures using statistical copulas: \cite{Garcia2016}; 
\item 
Testing (non)-informative sampling conditions in statistical surveys: \cite{Sikov2019}; 
\item 
Model selection for generalized Poisson distributions:
\cite{Hubert2009,SternZacks2002}; 
\item 
Model selection for generalized jump diffusion and Brownian motions, extremal distributions, and persistent memory processes: 
\cite{Andrade2015,Barahona2016,Bernardini2011,Kostrzewski2012,Rifo2009,Rifo2012}; 

\item 
Testing independence structures in contingency tables and multinomial models: \cite{Andrade2014,Bernardo2012,Oliveira2018,PereiraStern2008b}; 

\item 
Software certification according to compliance conditions: \cite{PereiraStern1999b}; 
\item 
Testing market equilibrium conditions for fundamental and financial derivative asset prices: \cite{Cerezetti2012}; 

\item 
Testing hypotheses in empirical economic studies: 
\cite{Chaiboonsri2018}; 
\item 
Event identification in acoustic signal processing: 
\cite{Hubert2018,Hubert2018b,Hubert2019};

\item 
Testing Hardy-Weinberg equilibrium in genetics: \cite{Brentani2011,Izbicki2012,Lauretto2009,Loschi2007,Montoya2001,Nakano2006,Rincon2010,Wittenburg2016}; 

\item 
Testing hypotheses in biological sciences, including cases in ecology, environmental sciences, medical diagnostics, efficacy evaluation of medical procedures, psychology and psychiatry: 
\cite{Ainsbury2013,Camargo2012,Kelter2020,Lima2014,Mathis2008,Pereira2a,Santos2020,Seixas2008,Shavitt2017,Spektor2019}; 

\item 
Testing hypotheses in astronomy and astrophysics: 
\cite{Chakrabarty2017,Johnson2009}. 

\item Comparison with p-values, other significance measures, and alternative hypothesis tests in miscellaneous applications: 
\cite{DCunha2016,Diniz2012b,KelterStern2020,Kelter2021a,Loschi2012,Patriota2013,Patriota2017}. 
  
\end{itemize}

\section{Epistemology, Ontology, and Metaphysics} 

Historically, the development of statistical significance measures (or logical truth-values), and the corresponding tests of hypotheses used in statistical science, have been reciprocally influenced by the epistemological frameworks in which they are presented. 
Rev.\,Thomas Bayes (1701-1761), whose work was communicated posthumously by Rev.\,Richard Price (1723-1791), developed the first methods in this area with clear goals in mind, as stated in \cite{Bayes1764}, also quoted in \cite[p.84]{Barker2001} and \cite[p.245]{Stern2018b}: 
\begin{quote} 
\textit{ 
The purpose, is to shew [show] what reason we have for 
believing that there are in the constitution of things, {fixed laws} according to which events happen, and that, therefore, the frame of the world must be the effect of the {wisdom and power of an intelligent cause}; and thus to confirm the argument taken from final causes for the existence of the Deity [...] \\  
It will be easy to see that the problem solved in this essay is more directly applicable to this purpose; for it shews [shows] us, with distinctness and precision, in every case of any particular order or recurrency of events, what reason there is to think that such recurrency or order is derived from {stable causes or regulations in nature}, and not from any of the irregularities of chance.} 
\end{quote} 

The next two generations in the development of probability and statistics, led by, among others, Pierre-Simon de Laplace (1749-1827) and George Boole (1815-1864), kept these core goals unchanged. 
While theological questions lost interest over time, the emphasis of statistical research remained essentially metaphysical -- in the (gnoseological) sense of searching for and justifying causal explanations for manifested phenomena. 
Moreover, these causal links were expressed as natural laws, ideally taking the form of mathematical equations, and later translated as sharp or precise hypotheses in statistical models, see   \cite{Stern2011a,Stern2011b,Stern2018b}. 
Karl Pearson (1857-1936) was the founder and leader of \textit{frequentist statistics}, the dominant school of thought in statistical science during the XX century. 
Influenced by the \textit{Inverted Spinozism} philosophy of Johann Gottlieb Fichte (1762-1814) and the \textit{Positivist} ideas of Auguste Comte (1798-1857), K.Pearson radically changed the goals of statistical science. 
He deprecated any form of causal reasoning, the conception or verification of natural laws, or the use of metaphysical (non-observable) entities. 
Instead, frequentist statistics, following a strict Positivist agenda, only aims to produce good-fitted empirical models able to describe or predict directly observed quantities. 
Accordingly, the very use of probability calculus is restricted to variables in the sample space, and strictly forbidden for (latent) variables in the parameter space, see \cite{Stern2018b}. 

Bruno de\,Finetti (1906–1985) is responsible for the resurgence of Bayesian statistics in the second half of the XX century, see \cite{Finetti1975,Finetti2006}. 
De\,Finetti reintroduced probability calculus in the parameter space. 
Latter on, this expanded use of probability language proved to be very useful for it could accommodate new means and methods provided by computer science, like Markov Chain Monte Carlo and other probabilistic algorithms, see \cite{Gilks1996,Hammersley1964,Ripley1987}. 
Nevertheless, philosophically, the de\,Finettian revolution was a very conservative one, remaining always amenable to the Positivist agenda. 
This was accomplished by making probabilities in the parameter space quantities of subjective and ephemeral character (integration variables), entities of low ontological status used only at intermediate steps in the computation of predictive probabilities for observable variables in the sample space, see \cite{Finetti1975,Finetti2006}. 

The \textit{Objective cognitive constructivism} epistemological framework was developed to host the e-value, the FBST, and their formal properties, including their ability to evaluate and test sharp statistical hypotheses. 
It also provides a naturalized way to ontology and metaphysics in empirical sciences via statistics, that is, a natural way to evaluate empirical support for natural laws and their accompanying causal explanations, and to validate the ``objective'' use of non-directly observable (or metaphysical) entities, in accordance with the goals of the original works of Bayes and Laplace. 
For further details on the Objective cognitive constructivism epistemological framework and the way in which, on the one hand, it accommodates the e-value, the FBST and their logical properties, and, on the other hand, it provides a naturalized approach to ontology and metaphysics, see \cite{Stern2007a,Stern2007b,Stern2008a,Stern2008b,Stern2011a,Stern2011b,Stern2014a,Stern2015,Stern2017a,Stern2017b,Stern2018a,Stern2018b,Stern2018c,Stern2020,Stern2022,Stern2023}.

\section{Future Research and Final Remarks}

In December 2018, at CLE-UNICAMP, Walter Carnielli organized the workshop \textit{Induction, Probability and their Dilemmas}, where the first author gave the presentation 
\textit{The Problem of Induction in Statistical Science}. 
During this workshop, we discussed several topics for further research, most of them motivated by areas of common interest that, nevertheless, are usually approached quite differently by the communities of Logic and Statistics. 
Four of these topics are presented in the sequel, followed by some additional topics for further research.

\subsubsection*{Statistics in (Un)Countable Sentential Probability} 

There is a long standing tradition in Logic to formalize probabilistic reasoning over finitary or countable sets of sentences in a language. 
This approach is appealing for its theoretical simplicity, for allowing computer-efficient implementations of specific models, etc. 
For historical analyzes of this approach, see \cite{Costa2001,Hailperin1996,Hailperin2011,Russell2015}, 
for recent works from Walter's group following this line, see \cite{Bueno2016,Rodrigues2021}. 
Notwithstanding their usefulness, simplicity, and popularity, finitary or countable sentential formalisms also have their limitations, for example, being unable to express measure-theoretic arguments used in mathematical statistics. 

Sections 6 and 7 sowed how the statistical definitions of the e-value and the GFBST characterize their logical properties, and hinted at how arguments of mathematical analysis enable one to travel the other way around. 
In this context, it is a topic for further research to explore how much of this theory and its applications can be expressed using the sentential probability approach, either by expanding the underlying languages (like infinitary, second order, or fixed-point logics), see \cite{Ebbinghaus1999,Gaifman1964,Gaifman1982,Hailperin1996,Hailperin2011,Hedman2004,Howson2009}, or by exploring relevant 
properties in specific statistical models (including topological properties, like countability and compactness, and regularity conditions of constraints and distributions, like bounded continuity and differentiability), see \cite{Borges2007,Dugundji1966,Finetti1975,Kelter2021b,Lipschutz1965,Searcoid2006}.

\subsubsection*{Functional Compositionality Structures} 

As already noticed in Section 6, the compositional structure of e-values can be characterized as a \textit{possibilistic abstract belief calculus}, as analyzed in 
\cite{Borges2007,Darwiche1992,Darwiche1993,Dubois1982,Dubois2012,Borges2007,Klir1988,Stern2003,SternPereira2014}. 
This particular algebraic formalism resembles in many ways the formalism used in statistical reliability theory for the analysis of complex systems assembled by serial/ parallel composition of simple(r) elements, see \cite{Barlow1981,Birnbaum1961,Kaufmann1977,Shier1991}.
As noticed in \cite[Sec.1]{Borges2007}, 
this formal resemblance was a source of inspiration at early stages of this research program. 
At the same time, the truth-function interpretation of the cumulative surprise function, $W(v)$, was inspired by the work of Ludwig Wittgenstein \cite{Cobde1998,Wittgenstein1921}. 
Finally, the concept of pragmatic hypotheses and its use in conjunction with the GFBST in testing scientific theories was inspired by sensitivity analysis, as used in optimization and systems' theory, and its interpretation in terms of fuzzy or paraconsistent logics, see \cite{Esteves2019,Stern2004}. 
While the many theoretical results already obtained in this research program justify (in our opinion) the analogies we made and the conceptual links we established between fragments of formal structures used in distinct (and often faraway) research areas, we believe that these intuitive connections deserve and can benefit from an even more rigorous and general setting. 
In this context, the tools of category theory or other formal abstraction methods, see \cite{Carnielli2008a,Carnielli2008c,Fong2019,Lawvere1997,Spivak2014}, offer an opportunity for further research in the study of basic logical properties of the aforementioned systems, specially concerning the investigation of their essential compositionality structures and possible generalizations. 


\subsubsection*{Rough and Fuzzy Sets} 
Essential logical properties of the GFBST can be explained regarding the e-value as a transformation between probability and possibility measures, see \cite{SternPereira2014}. 
The underlying interconnections between alternative representations of uncertainty at the core of this theory 
provide a general motivation to more specific topics of further research presented in this section. 

In the GFBST framework, a sharp hypothesis can be either rejected or remain undecided, but can never be accepted, see Section 7.
Aiming to allow the acceptance of surrogate hypotheses of interest, the authors showed how to enlarge an underlying sharp hypothesis into a slack \textit{pragmatic hypothesis}, see \cite{Esteves2019}. 
This is accomplished by taking into account (im)precisions of measurement equipment, (un)certainties of fundamental constants, and other relevant metrological or methodological error bounds. 
It should be remarked that the pragmatic hypotheses defined in \cite{Esteves2019} are crisp sets. 

In the context of exact sciences, the sharp (and, hence, crisp) nature of the original statistical hypothesis of interest is well-supported by the \textit{Objective cognitive constructivism} epistemological framework; see \cite{Stern2007a,Stern2007b,Stern2008a,Stern2011a,Stern2011b,Stern2014a,Stern2015,Stern2017a,Stern2017b,Stern2018b,Stern2018c,Stern2020}. 
In contrast, the definition of the associated pragmatic hypothesis as a crisp set may be considered an over-simplification. 
Rough and fuzzy sets offer a theoretical framework that allows the construction of generalized pragmatic hypotheses, see 
\cite{Dubois1982,Dubois2004,Dubois2007,Dubois2012,Klir1988,Klir1995,Salicone2007}. 
Rough or fuzzy pragmatic hypotheses should be able to overcome artificial limitations imposed by crisp set representations, support natural and intuitive interpretations in varied contexts of application, and still preserve the best logical properties of the GFBST. 

The theoretical framework of paraconsistent logic can be used to interpret sensitivity analyses of the e-value concerning changes in the statistical model's prior or reference measures, see \cite{Stern2004}. 
The same article shows how to integrate such sensitivity tests into (crisp) confidence intervals. 
Following the same rationale used in the last paragraph, these crisp intervals could be generalized to rough or fuzzy sets able to provide better representations of pertinent uncertainties.

\subsubsection*{Law, Complexity, and (In)Consequence in Social Systems} 
Early scientists used the expression \textit{natural law} as a metaphor, that ferried to the domain of nature the normative character that a law has over human behavior. 
Now that mankind learned a good deal about natural laws, including their mathematical, logical, epistemic, ontological, and metaphysical characteristics, we could perhaps travel the metaphorical path taken by early scientists in the opposite direction. 
We could do so in an attempt to use lessons learned about natural laws to better understand the nature of laws intended to establish norms for human behavior, relations, and interactions in social systems. 
We believe this task can be accomplished within the framework provided by Niklas Luhmann (1927–1998) \textit{Sociological theory of law}, see \cite{Campilongo2020,Foerster2001,Foerster2003,Luhmann1985,Luhmann1989,Maturana1980,Rasch2000,Stern2018a,Teubner1988,Zeleny1979}. 
  
In Luhmann's theory of law, the main purpose of the legal system is \textit{congruent generalization of normative behavior expectations}, see \cite[pp.77,82]{Luhmann1985}. 
This emblematic statement can be interpreted as follows, see \cite{Stern2018a}: 
In Luhmann's view, norms are neither preexisting conditions nor a priori factual realities. 
Instead, norms are conceived as intentional projects or idealized models of how society should be; see \cite[p.40]{Luhmann1985}. 
Moreover, in such idealized models, social harmony is based on establishing well-defined, stable, sustainable, and reliable behavioral patterns, also known as \textit{eigen-behaviors} or \textit{eigen-solutions}, see \cite{Foerster2001,Stern2018a}. 
Furthermore, laws and regulations of a society should reflect its norms, stimulating/ penalizing forms of conduct that sustain/ disrupt virtuous eigen-behaviors. 
Finally, in Luhmann's view, social norms are not static but essentially dynamic, co-evolving (over long runs) with the behavioral patterns in the society that they simultaneously try to describe, regulate and stabilize (in their present forms, over short runs). 
In this context, we can try to reinterpret some lessons about model selection and their complexity learned in Section 9, claiming that, as in the case of good empirical laws, good social laws should follow the golden path of equilibrium, avoiding extremes of scarcity and excess, see also \cite{Stern2007b,Stern2011b}. 
  
On the one hand, oversimplified social laws fail to capture important distinctions considered necessary or relevant for establishing sustainable eigen-behaviors. 
On the other hand, excessively complex legislation creates all sorts of misinterpretations, unforeseen loopholes, and other unintended consequences. 
Moreover, such spurious side effects may not only obstruct virtuous eigen-behaviors, but even induce vicious ones (that must then be detected, identified, and inhibited). 
Furthermore, increasing legal complexities imply increasing processing times, delayed justice, and greater economic costs to operate the legal system -- all burdens to be paid by the society the same system serves. 
  
The areas of Logic and Computer science have developed several methods to measure computational complexity, either by counting processing operations in an algorithm or by accessing the code length of its description, see 
\cite{Ebbinghaus1999,Epstein2008,Inhasz2010,Rissanen1989,Wallace1968,Wallace1999}. 
Some of the methods used to measure complexity in statistical models have already been explained in Section 9. 
Meanwhile, Information science and Systems' theory have focused on entropy related measures of complexity, see 
\cite{Aczel1975,Attneave1959,Dugdale1996,Kapur1989,Stern2011a,Stern2014b,Takada2020,Wallace2005}. 
In contrast, the mathematical treatment of complexity in theoretical and empirical legal studies is still incipient, see 
\cite{Cauble2021,Kades1997,Muchmore2016,Murray2019,Ruhl2008,Ruhl2015,Sichelman2021,Schuck1992,Sichelman2021}. 
Na\"{\i}ve adaptations of complexity measures artificially borrowed from other areas often fail to capture relevant aspects of the legal environment. 
  
Although this topic of further research can potentially benefit from all the sophisticated technical developments and approaches to complexity theory already mentioned, 
it also requires a healthy dose of critical thinking applied to the human condition in daily life, yet another area of interest of Walter's research group, see \cite{Carnielli2010}.

\subsubsection*{Individual vs. Collective Liability, Legal Burden of Proof, etc.}

As commented at the end of Section 2, the \textit{Onus Probandi} and \textit{In Dubio Pro Reo} principles, as they are used in the legal system, were a source of inspiration for the definitions of the e-value and the FBST, see \cite{Gaskins1992,Guy1994,Kokott1998,PereiraStern1999b,PereiraNS2000,Pigliucci2014,Stern2003,Stern2018a}. 
Moreover, it has long been recognized that many decision procedures traditionally used in Bayesian statistics are incompatible with these principles. 
Dennis Lindley uses the famous \textit{gatecrasher} example to argue that, under certain circumstances, collective responsibility is a preferable principle of justice, and that legal decision procedures implying collective liabilities and punishments should be accepted, 
see \cite{Lindley1991,Cohen1977}. 
Furthermore, Lindley shows how legal decisions based on Bayes Factors are compatible with such a collective liability principle.  
   
The right tool for the right job! 
In accordance with this dictum, further research could contrast the legal and social implications of alternative liability principles, their potential or desirability, and their compatibility with alternative inference theories, decision procedures, and epistemological frameworks.

\subsubsection*{e-value and FBST under Non-Standard Regularity Conditions} 

The asymptotic consistency and standardization procedures presented in Section 5 rely on standard regularity conditions that include continuity and differentiability of several functions defining the statistical model, and also the location of the true value of the parameter at the topological interior of the parameter space. 
However, many interesting models present non-standard conditions. 
For example, mixture models and separate hypotheses models often present hypotheses at the (topological) border of the parameter space, 
see \cite{Lauretto2005a,Lauretto2007}, whereas eigenvalue and eigenvector inference problems may present far more challenging conditions, see \cite{Schwartzman2008}. 
The methods presented in \cite{Andrews2001,Chernoff1954,Drton2009,Self1987} 
and references therein provide tools to generalize the e-value standardization procedure and the corresponding asymptotic results to non-standard regularity conditions.     

Other non-standard inference conditions arise in semi-parametric models like, for example, Hilbert space models based on Fourier, orthogonal polynomials, wavelets, and other infinite functional bases. 
In these cases, it is often necessary to use weakly informative priors that dampen high-frequency modes and that, in this way, effectively shrink resulting truncated models, see \cite[Sec.6]{Stern2020} and references therein. 
Further research should explore proper theory and convenient means and methods needed to use the e-value and the FBST under a variety of non-standard regularity conditions.

\subsubsection*{Efficient Computational Methods and Implementations} 

The e-value and FBST research program would greatly benefit from computational tools tailor-made for computational efficiency, reliability, and ease of use in the numerical tasks described at Section 8. 
Most beneficial would be the construction of an open-source and user-friendly environment combining specially adapted versions of the following techniques:    
(1) Efficient Monte-Carlo integrators based on Hit-and-Run, Nested Sampling and similar algorithms, 
see \cite{Chen2000,Karawatzki2005,Skilling2004,Skilling2006};  
(2) Efficient computational condensation procedures used in conjunction with the aforementioned Monte-Carlo procedures,  
see \cite{Kaplan1987,Williamson1989};   
(3) Efficient stochastic optimization methods working in conjunction with the aforementioned Monte-Carlo procedures, see  \cite{Pflug1996,Spall2003,Stern1992,Voigt1990,Wah2008};  
(4) Higher order asymptotic approximations used either to develop stand-alone fast e-value calculation procedures or to develop variance reduction techniques for the aforementioned Monte-Carlo procedures, see \cite{Cabras2015,Pinto2012,Ranzato2018,Ruli2020,Ruli2020b,Ventura2013,Ventura2014,Ventura2016};   
(5) Efficient numerical convolution procedures, used in conjunction with the aforementioned condensation procedures, as required in the analysis of complex models as described in Section 6.

\subsubsection*{Acknowledgments} 

The authors take this opportunity to thank Prof. Walter Carnielli for, along the years, inviting, welcoming, and stimulating investigations, discussions, and interactions at areas of overlap between logic, probability, statistics, formal methods in science, philosophy of science, and general philosophy. 
The authors have benefited greatly from Walter Carnielli's encyclopedic knowledge, open-mindedness, positive attitude, and active engagement in fostering research in all these areas, either at his home base at CLE - the Center for Logic, Epistemology, and History of Science of UNICAMP, or throughout Brazil and around the world.  
The authors are grateful to Profs. Abilio Rodrigues,  
Alfredo Freire and Henrique Antunes 
for organizing a Festschrift celebrating Walter's 70th birthday. 
 
The authors thank the organizers of the following conferences, where this paper was presented:  
CLMPST-2023, the 17th Congress on Logic, Methodology and Philosophy of Science and Technology, at FCE-UBA, Buenos Aires, July 24-29, 2023;  
Principia-2023, the 13th Principia International Symposium, at UFSC, Florian\'{o}polis, August 14-17, 2023; 
RatioLog-2024, the 2nd Workshop on Logic, Rationality, and Probability, at CLE-Unicamp, Campinas, April 22-25, 2024; and   
MaxEnt-2024, the 43rd International Workshop on Bayesian Inference and Maximum Entropy Methods in Science and Engineering, at U-Ghent, July 1-5, 2024.

\subsubsection*{Support, Funding, and Conflicts of Interest}

The authors have received direct or indirect support from the following institutions: 
USP -- University of S\~{a}o Paulo (JMS, CABP, MSL, LGE, and WSB); 
UFSCar -- Federal Univesity of S\~{a}o Carlos (RI, RBS); 
ABJ -- Brazilian Jurimetrics Association (JMS, RBS); 
CNPq -- Brazilian National Counsel of Technological and Scientific Development, grants 
PQ 303290/2021-8 (JMS),
PQ 302767/2017-7 (CABP), and 
PQ 309607-2020-5 (RI); 
FAPESP -- State of S\~{a}o Paulo Research Foundation, grants 
2019/11321-9 (RI), 
CEPID CeMEAI -- Center for Mathematical Sciences Applied to Industry 2013/07375-0 (JMS, MSL), and 
CEPID Shell-RCGI -- Research Center for Greenhouse Gas Innovation 2014/50279-4 (JMS); and 
CAPES -- Coordination for Improvement of Higher Education Personnel, finance code 001 (RI, RBS, MAD). 
Concerning this article, the authors have no conflicts of interest to declare. 

\mbox{} \vspace{2mm} \mbox{} 


\setlength{\columnsep}{15pt} 
\begin{multicols}{2} 


\def\SBlongestlabel{}
\SBtitlestyle{simple}
\SBsubtitlestyle{simple} 
       

 \begin{small} 


 \end{small} 
        
\end{multicols}


\begin{thebibliography}{}          
     
         
\bibitem{Aczel1975}
Acz\'{e}l, J\'{a}nos Dezs\"{o}; Dar\'{o}czy, Zolt\'{a}n. 
\textit{On Measures of Information and Their Characterizations}. 
NY, Academic Press, 1975. 
       
\bibitem{Ainsbury2013} 
Ainsbury, Elizabeth A; Vinnikov, Volodymyr A; Puig, Pedro; Higueras, Manuel; Maznyk, Nataliya A; Lloyd, David C;  Rothkamm, Kai. 
Review of Bayesian statistical analysis methods for cytogenetic radiation biodosimetry, with a practical example. 
\textit{Radiation Protection Dosimetry}, 162, 3, 185-196, 2013.  

\bibitem{Amari1987} 
Amari, Shun Ichi; Barndorff-Nielsen, Ole Eiler; Kass, Robert E; Lauritzen, Steffen L; Rao, Calyampudi Radhakrishna; Gupta, Shanti S.  
\textit{Differential Geometry in Statistical Inference}. 
Hayward, IMS, 1987. 
   
\bibitem{Amari2007} 
Amari, Shun Ichi; Nagaoka, Hiroshi. \textit{Methods of Information Geometry}. Providence, RI, AMS, 2007. 

\bibitem{Amari2016} 
Amari, Shun Ichi. 
\textit{Information Geometry and Its Applications}. 
Berlin, Springer, 2016. 

\bibitem{Andrade2015}
Andrade, Plinio; Rifo, Laura Leticia Ramos; Torres, Soledad;   Torres-Avil\'{e}s, Francisco. 
Bayesian Inference on the Memory Parameter for Gamma-Modulated Regression Models. 
\textit{Entropy}, 17, 10,  6576-6597, 2015.  

\bibitem{Andrade2014}
Andrade, Pablo De Morais; Stern, Julio Michael; Pereira, Carlos Alberto De Bragan\c{c}a.  
Bayesian Test of Significance for Conditional Independence: The Multinomial Model. 
\textit{Entropy}, 16, 3, 1376-1395, 2014.   

\bibitem{Andrews2001} 
Andrews, Donald W. K. 
Testing when a Parameter in not on the Boundary of the Maintened Hypothesis. 
\textit{Econometrica}, 69, 3, 683-734, 2001.  


\bibitem{Assane2018}   
Assane, Cachimo Combo; Pereira, Basilio de Bragan\c{c}a; Pereira, Carlos Alberto de Bragan\c{c}a.  
Bayesian significance test for discriminating between survival distributions. 
\textit{Communications in Statistics - Theory and Methods}, 47, 24, 6095-6107, 2018. 

\bibitem{Assane2019}   
Assane, Cachimo Combo; Pereira, Basilio de Bragan\c{c}a; Pereira, Carlos Alberto de Bragan\c{c}a.  Model choice in separate families: A comparison between the FBST and the Cox test. 
\textit{Communications in Statistics - Simulation and Computation}, 48, 9, 2641-2654, 2019. 

\bibitem{Attneave1959}  
Attneave, Fred. Applications of Information Theory to Psychology: A summary of basic concepts, methods, and results. NY, Holt, Rinehart and Winston, 1959. 


\bibitem{Barker2001} 
Barker, Frank; Evans Robin.  
\textit{The Probability of Mr Bayes: A constructive re-evaluation of Mr Bayes' essay and of the opinions concerning it expressed by various authorities}. Technical report, Department of Electrical and Electronic Engineering, University of Melbourne, 2001.  

\bibitem{Barlow1981}
Barlow,  Richard E; Prochan, Frank. 
\textit{Statistical Theory of Reliability and Life Testing Probability Models}. Silver Spring, To Begin With, 1981. 

   
\bibitem{Barahona2016} 
Barahona, Manuel; Rifo, Laura; Sep\'{u}lveda, Maritza; Torres, Soledad.   
A Simulation-Based Study on Bayesian Estimators for
the Skew Brownian Motion.  
\textit{Entropy}, 18, 7, 241, 1-14, 2016.   
  
\bibitem{Barak2012}  
Barak, Aharon. \textit{Proportionality: Constitutional Rights and Their Limitations}. 	Cambridge Univ. Press, 2012.  
     
\bibitem{Basu1988}  
Basu, Debabrata; Ghosh, J.K.  
\textit{Statistical Information and Likelihood}. 
Lecture Notes in Statistics, 45, 1988. 

\bibitem{Bayes1764}
Bayes, Thomas. An Essay Towards Solving a 
Problem in the Doctrine of Chances. \textit{Phil. Trans. Roy. Soc. London}, 53, 370-418, 1764. 


\bibitem{Berger1988}
Berger, James O; Wolpert, Rober L.  \textit{The Likelihood Principle}, 2nd ed. Hayward, CA, Inst of
Mathematical Statistic, 1988. 

\bibitem{Bernardini2011}
Bernardini, Diego F. de; Rifo, Laura Leticia Ramos. Full Bayesian significance test for extremal distributions. 
\textit{Journal of Applied Statistics}, 38, 4, 851-863, 2011. 

\bibitem{Bernardo2012} 
Bernardo, Gustavo; Lauretto, Marcelo de Souza; Stern, Julio Michael. The Full Bayesian Significance Test for Symmetry in Contingency Tables. 
\textit{AIP Conference Proceedings}, 1443, 198-205, 2012.

\bibitem{Bernardo2005}  
Bernardo, Jos\'{e} M.  
Reference Analysis. p.17-90 in  
Dey, D.K; Rao, C.R. Bayesian Thinking: 
Modeling and Computation. 
\textit{Handbook of Statistics}, 
v.25. Amsterdam, Elsevier, 2005.   

\bibitem{Beziau2012}
B\'{e}ziau, Jean-Yves. 
The power of the hexagon. \textit{Logica Universalis}, 6, 1-43, 2012. 

\bibitem{Birgin2014}
Birgin, Ernesto Julian Goldberg; Mart\'{\i}nez, Jos\'{e} Mario. \textit{Practical Augmented Lagrangian Methods for Constrained Optimization}. SIAM, 2014. 

\bibitem{Birnbaum1961} 
Birnbaum, Zygmunt Wilhelm; Esary, James D; Saunders, Samuel C.  Multicomponent Systems and Structures,
and their Reliability. \textit{Technometrics}, 3, 55-77, 1961. 

\bibitem{Blanche1966}
Blanch\'{e}, Robert. \textit{Structures Intellectuelles: Essai sur l'Organisation Syst\'{e}matique des Concepts}. Vrin, 1966. 

 
\bibitem{Borges1979} 
Borges, Wagner de Souza. \textit{Modelos Probabil\'{\i}sticos de Confiabilidade}. Rio de Janeiro, IMPA -- Instituto de Matem\'{a}tica Pura e Aplicada, 1979. 
  
   
\bibitem{Borges2006} 
Borges, Wagner de Souza; Stern, Julio Michael. Evidence and Compositionality. pp.307-315 in Lawry, J.\,et\,al. 
\textit{Soft Methods for Integrated Uncertainty Modelling}, 
Berlin: Springer. Proceedings of SMPS-2006. 
 
\bibitem{Borges2007} 
Borges, Wagner de Souza; Stern, Julio Michael. The Rules of Logic Composition for the Bayesian Epistemic E-Values. 
\textit{Logic Journal of the IGPL}, 15, 5/6, 401-420, 2007.

\bibitem{Box1973} 
Box, George Edward Pelham; and Tiao, George C.  
\textit{Bayesian Inference in Statistical Analysis}. 
London, Addison-Wesley, 1973.

\bibitem{Bracewell1986} 
Bracewell, Ronald N. \textit{The Fourier Transform and its Aplications}. McGraw-Hill, 1986. 

\bibitem{Brentani2011}
Brentani, Helena; Nakano, Eduardo Y; Martins, Camila B; Izbicki, Rafael; Pereira, Carlos Alberto. 
Disequilibrium Coefficient: A Bayesian Perspective. 
\textit{Statistical Applications in Genetics and Molecular Biology}, 10, 1, 22, 1-24, 2011. 

\bibitem{Brooks2011} 
Brooks, Steve; Gelman, Andrew; Jones, Galin L; Meng, Xiao-Li. \textit{Handbook of Markov Chain Monte Carlo}. CRC Press, 2011. 

\bibitem{Bueno2016} 
Bueno-Soler, Juliana; Carnielli, Walter Alexandre.
Paraconsistent Probabilities: Consistency, Contradictions and Bayes' Theorem. \textit{Entropy}, 18, 325, 1-18, 2016.

\bibitem{Cabras2015} 
Cabras, Stefano; Racugno, Walter; Ventura, Laura.  
Higher order asymptotic computation of Bayesian significance tests for precise null hypotheses in the presence of nuisance parameters. 
\textit{Journal of Statistical Computation and Simulation}, 85, 15, 2989-3001, 2015.     

\bibitem{Campilongo2020}
Campilongo, Celso Fernandes; Amato, Lucas Fucci;  
Barros, Marco Antonio Loschiavo Leme de. 
\textit{Luhmann and Socio-Legal Research: An Empirical Agenda for Social Systems}. NY, Routledge, 2020.  


\bibitem{Camargo2012} 
Camargo, Andr\'{e} P; Stern, Julio M; Lauretto, Marcelo S. Estimation and model selection in Dirichlet regression. 
\textit{AIP Conference Proceedings}, 1443, 206-213, 2012. 
   
\bibitem{Carnielli2008a}
Carnielli, Walter Alexandre; Pizzi, Claudio. 
\textit{Modalities and Multimodalities}. 
Springer-Science, 2008. 

\bibitem{Carnielli2008c} 
Carnielli, Walter Alexandre; Coniglio, Marcelo Esteban; Gabbay, Dov; Gouveia, P; Sernadas, C. \textit{Analysis and Synthesis of Logics - How To Cut And Paste Reasoning Systems}. NY, Springer, 2008. 
  
\bibitem{Carnielli2010} 
Carnielli Walter Alexandre; Epstein, Richard L.  
\textit{Pensamento Cr\'{\i}tico}. 
S\~{a}o Paulo, Rideel, 2010. 

\bibitem{Cauble2021}
Cauble, Emily. Unsophisticated Taxpayers, Rules Versus Standards, and Form Versus Substance. \textit{Loyola University Chicago Law Journal}, 52, 329-349, 2021. 

\bibitem{Cerezetti2012} 
Cerezetti, Fernando Valvano; Stern, Julio Michael. Non-arbitrage in Financial Markets: A Bayesian Approach for Verification. 
\textit{AIP Conference Proceedings}, 1490, 87-96,  2012. 

\bibitem{Chakrabarty2017} 
Chakrabarty, Dalia. A New Bayesian Test to Test for the Intractability-Countering Hypothesis. 
\textit{Journal of the American Statistical Association}, 112, 518, 561-577, 2017.  

\bibitem{Chaiboonsri2018}  
Chaiboonsri, Chukiat;  Wannapan, Satawat; Saosaovaphak, Anuphak.  
Economic and Business Cycle of India:
Evidence from ICT Sector. 
p.29-43 in Tsounis, Nicholas; Vlachvei, Aspasia. 
\textit{Advances in Panel Data Analysis in Applied Economic Research}. 
Cham, Switzerland, Springer Nature, 2018.     

\bibitem{Chen2015} 
Chen, C. W. S; Lee, S. A local unit root test in mean for financial time series. 
\textit{Journal of Statistical Computation and Simulation}, 86, 4, 788-806, 2015.

\bibitem{Chen2000}  
Chen, Ming-Hui; Shao, Qi-Man; Ibrahim, Joseph G.	
\textit{Monte Carlo Methods in Bayesian Computation}. 
Springer-Verlag, 2000. 

\bibitem{Cherkassky1998}
Cherkassky, Vladimir; Mulier, Filip M. \textit{Learning from Data}. NY, Wiley, 1998.  

\bibitem{Chernoff1954}
Chernoff, Herman. On the Distribution of the Likelihood Ratio. 
\textit{Ann.\,Math.\,Statist.} 25, 3, 573-578, 1954.  


\bibitem{Cohen1977} 
Cohen, Jonathan.  
\textit{The Probable and the Provable}. Oxford: Clarendon, 1977.  

\bibitem{Cobde1998}
Cond\'{e}, Mauro L\'{u}cio Leit\~{a}o. \textit{Wittgenstein, Linguagem e Mundo}. S\~{a}o Paulo, Annablume, 1998. 

\bibitem{Costa2001} 
Costa, Newton Carneiro da; Tsuji, Marcelo. Review of Theodore Hailperin, Sentential Probability Logic.  
\textit{Modern Logic}, 8, 3-4, 103-107, 2001. 
  


\bibitem{Darwiche1992}
Darwiche, Adnan Youssef;  Ginsberg, Matthew L.  A Symbolic Generalization of Probability Theory. \textit{AAAI-92. 10-th Conf. American Association for Artificial Intelligence}, 1992. 

\bibitem{Darwiche1993}
Darwiche, Adnan Youssef. \textit{A Symbolic Generalization of Probability Theory}. Ph.D. Thesis,
Stanford University, 1993.

\bibitem{DCunha2016} 
D'Cunha, Juliet Gratia; Rao, Aruna  K.  
Frequentist Comparison of the Bayesian Significance Test for Testing the Median of the Lognormal Distribution. 
\textit{InterStat}, 2016, 02, 001, 1-25, 2016.  

\bibitem{DeGroot1970} 
DeGroot, Morris Herman. Optimal Statistical Decisions. 
NY, McGraw-Hill, 1970. 

\bibitem{DeGroot2011} 
DeGroot, Morris H; Schervish, Mark J.  
\textit{Probability and Statistics}. Pearson, 2011.  

\bibitem{Diniz2011} 
Diniz, Marcio; Pereira, Carlos Alberto de Bragan\c{c}a; Stern, Julio Michael. Unit Roots: Bayesian Significance Test. 
\textit{Communications in Statistics - Theory and Methods},  40, 23, 4200-4213, 2011. 

\bibitem{Diniz2012a} 
Diniz, Marcio; Pereira, Carlos Alberto de Bragan\c{c}a; Stern, Julio Michael. 
Cointegration: Bayesian Significance Test. 
\textit{Communications in Statistics - Theory and Methods}, 41, 19, 3562-3574, 2012. 

\bibitem{Diniz2012b} 
Diniz, Marcio; Pereira, Carlos A.B; Polpo, Adriano; 
Stern, Julio M; Wechsler, Sergio. 
Relationship between Bayesian and Frequentist Significance Indices. 
\textit{International Journal for Uncertainty Quantification}, 2, 2, 161-172, 2012.   

\bibitem{Diniz2022} 
Diniz, Marcio; Pereira, Carlos Alberto de Bragan\c{c}a; Stern, Julio Michael. 
\textit{e-value} entry at Wiley's \textit{StatsRef: Statistics Reference}, 2022.  
\texttt{doi:10.1002/9781118445112.stat08375} 

\bibitem{Drton2009}
Drton, Mathias. Likelihood Ratio Tests and Singularities. 
\textit{Ann.\,Statist.}, 37, 2, 979-1012, 2009. 
        

\bibitem{Dubois1982}   
Dubois, Didier; Prade, Henri. On Several Representations of an Uncertain Body of Evidence. 
pp.167-181 in Gupta, M.M;  Sanchez, E; eds.  
\textit{Fuzzy Information and Decision Processes}, 
North-Holland, 1982.

\bibitem{Dubois2004}
Dubois,Didier; Foulloy, Laurent; Mauris, Gilles; Prade, Henri. Probability-possibility transformations, triangular fuzzy sets, and probabilistic inequalities. \textit{Reliable Computing}, 10, 273–297, 2004. 

\bibitem{Dubois2007} 
Dubois, Didier; Prade, Henri. 
Rough Fuzzy Sets and Fuzzy Rough Sets. 
\textit{International Journal of General Systems}, 2007.  
\texttt{doi:10.1080/03081079008935107}   

\bibitem{Dubois2012}
Dudois, Didier; Prade, Henri. 
From Blanch\'{e}'s Hexagonal Organization of Concepts to Formal Concept Analysis and Possibility Theory. 
\textit{Logica Universalis}, 6, 149-169, 2012. 
   
\bibitem{Dugdale1996} 
Dugdale, J. Sydney. Entropy and Its Physical Meaning. London, Taylor and Francis, 1996. 
 
\bibitem{Dugundji1966}
Dugundji, James. \textit{Topology}. Boston, Allyn and Bacon, 1966.  

\bibitem{Dumitriu1977} 
Dumitriu, Anton. \textit{History of Logic}. Abacus Press, 1977. 

\bibitem{Eaton1989} 
Eaton, Morris L. \textit{Group Invariance Applications in Statistics}. Hayward, CA, IMS, 1989. 

\bibitem{Ebbinghaus1999} 
Ebbinghaus, Heinz-Dieter; Flum, J\"{o}rg. \textit{Finite Model Theory}.  Berlin, Springer, 1999.

\bibitem{Engle2012} 
Engle, Eric Allen. The History of the General Principle of Proportionality. \textit{Dartmouth Law Journal}, 10, 1, 1–11, 2012. 

\bibitem{Epstein1974} 
Epstein, George; Frieder, Guideon; Rine, David. The Development of Multiple-Valued Logic as Related to Computer Science. \textit{Computer}, 7, 9, 20-32, 1974. 
  
\bibitem{Epstein1993}  
Epstein, George. \textit{Multiple-Valued Logic Design}. IOP Publishing, 1993. 
   
\bibitem{Epstein2008} 
Epstein, Richard L; Carnielli Walter Alexandre. 
\textit{Computability: Computable Functions, Logic, and the Foundations of Mathematics}. 
Socorro, NM, Advanced Reasoning Forum, 2008. 
      
\bibitem{Esteves2016} 
Esteves, Luis Gustavo; Izbicki, Rafael; Stern, Julio Michael; Stern, Rafael Bassi. 
The logical consistency of simultaneous agnostic hypothesis tests. 
\textit{Entropy}, 18, 256, 2016 

\bibitem{Esteves2019} 
Esteves, Luis Gustavo; Izbicki, Rafael; Stern, Julio Michael; Stern, Rafael Bassi. 
Pragmatic Hypotheses in the Evolution of Science. 
\textit{Entropy}, 21, 9, 883, 2019.  


\bibitem{Esteves2023} 
Esteves, Luis Gustavo ; Izbicki, Rafael; Stern, Julio Michael; Stern, Rafael Bassi. Logical Coherence in Bayesian Simultaneous Three-Way Hypothesis Tests. \textit{International Journal of Approximate Reasononig}, 152, 297-309, 2023. 

\bibitem{Evans1997} 
Evans, Michael. Bayesian Inference Procedures Derived via the Concept of Relative Surprise. 
\textit{Communications in Statistics}, 26, 1125-1143, 1997.


\bibitem{Fang1997}  
Fang, Shu Cheng; Rajasekera, Jay R; Tsao, Jacob, H.S.    
\textit{Entropy Optimization and Mathematical Programming}.
Kluwer, Dordrecht, 1997. 

\bibitem{Finetti1957} 
Finetti, Bruno de. \textit{Matematica Logico-Intuitiva}. Roma, Cremonese, 1957. 

\bibitem{Finetti1975}
Finetti, Bruno de. 
\textit{Theory of Probability}. New York, Wiley, 1997.

\bibitem{Finetti2006}
Finetti, Bruno de. \textit{L'invenzione della verit\`{a}}.  Milano, Scienza e Idee, 2006. 

\bibitem{Flanders1974}   
Flanders, Harley; Korfhage, RobertR; Price, Justin J. \textit{A Second Curse in Calculus}. NY, Academic Press, 1974.  

\bibitem{Foerster2003} 
- Foerster, Heinz von. \textit{Understanding Understanding: Essays on Cybernetics and Cognition}. NY, Springer Verlag, 2003.

\bibitem{Foerster2001}   
Foerster, Heinz von; Segal, Lynn. The Dream of Reality. Heintz von Foerster's Constructivism. NY, Springer, 2001. 

\bibitem{Fong2019} 
Fong, Brendan; Spivak, David Isaac. \textit{An Invitation to Applied Category Theory: Seven Sketches in Compositionality}. Cambridge Univ. Press, 2019
   

\bibitem{Fossaluza2017} 
Fossaluza, Victor; Izbicki, Rafael; Silva, Gustavo Miranda da;  Esteves, Lu\'{\i}s Gustavo.  
Coherent Hypothesis Testing. 
\textit{The American Statistician},  71, 3, 242-248, 2017. 


\bibitem{Frieden2004} 
Frieden, Roy Bernard \textit{Science from Fisher Information: A Unification}. Cambridge Univ. Press, 2004. 

\bibitem{Gaifman1964}
Gaifman, Haim. Concerning Measures in First Order Calculi. \textit{Israel J.\,of Mathematics}, 2, 1–18, 1964. 

\bibitem{Gaifman1982}
Gaifman, Haim; Snir, Marc. Probabilities Over Rich Languages, Testing and Randomness. 
\textit{The Journal of Symbolic Logic}, 47, 3, 495-548, 1982. 

\bibitem{Gallais1982} 
Gallais, Pierre.   \textit{Dialectique du R\'{e}cit Medi\'{e}val: Chr\'{e}tien de Troyes et l’Hexagone Logique}. Amsterdam, Rodopi, 1982. 



\bibitem{Garcia2016} 
Garc\'{\i}a, Jes\'{u}s E; Gonz\'{a}lez-L\'{o}pez, Ver\'{o}nica; Nelsen, Roger B.  
The Structure of the Class of Maximum Tsallis-Havrda-Chavat Entropy Copulas. 
\textit{Entropy}, 18, 7, 264, 1-6, 2016.   

\bibitem{Garson2006}
Garson, James W. \textit{Modal Logic for Philosophers}. Cambridge Univ. Press, 2006. 

\bibitem{Gaskins1992}
Gaskins, Richard H. Burdens of Proof in Modern Discourse. New Haven, Yale Univ. Press, 1992.

\bibitem{Gelman2004}
Gelman, Andrew; Carlin, John B; Stern, Hal S; Rubin, Donald B. 
\textit{Bayesian Data Analysis}. Chapman-Hall/ CRC, 2004.  

\bibitem{Gilks1996}
Gilks, Wally R; Richardson, Sylvia; Spiegelhalter, David John. \textit{Markov Chain Monte Carlo in Practice}. CRC Press, 1996. 

\bibitem{Good1983}
Good, Irving John. \textit{Good Thinking}. Univ. of Minnesota. 1983. 

\bibitem{Gottenwald2001} 
Gottenwald, Siegfried. \textit{A Treatise on Many-Valued Logics}. Research Studies Press, 2001. 

\bibitem{Guy1994}
Guy, Dan M; Carmichael, Douglas; Whittington, O. Ray. \textit{Audit Sampling: An Introduction}. John Wiley, 1994.  

\bibitem{Hailperin1996}
Hailperin, Theodore. \textit{Sentential Probability Logic: Origins, Development, Current Status, and Technical Applications}. Bethlehem, Lehigh Univ. Press, 1996.  

\bibitem{Hailperin2011}
Hailperin, Theodore. \textit{Logic with a Probability Semantics: Including Solutions to Some Philosophical Problems}. Bethlehem, Lehigh Univ. Press, 2011.  
  
\bibitem{Hammersley1964}
Hammersley, John Michael; Handscomb, David Christopher. \textit{Monte Carlo Methods}. Chapman and Hall, 1964. 

\bibitem{Hedman2004} 
Hedman, Shawn. \textit{A First Course in Logic: An Introduction to Model Theory, Proof Theory, Computability, and Complexity}. Oxford Univ. Press, 2004. 

\bibitem{Hocking1985} 
Hocking, Ronald R. The Analysis of Linear Models. Monterey: Brooks Cole, 1985.  

\bibitem{Howson2009}
Howson, Colin. Can Logic be Combined with Probability? Probably. \textit{Journal of Applied Logic}, 7, 2, 177–187, 2009. 
 
\bibitem{Hubert2009}  
Hubert, Paulo;  Lauretto, Marcelo de Souza;  Stern, Julio Michael 2009. 
FBST for Generalized Poisson Distribution. 
\textit{AIP Conference Proceedings}, 1193, 210-2019, 2009. 

\bibitem{Hubert2018}  
Hubert, Paulo; Padovese, Linilson;  Stern, Julio Michael. A Sequential Algorithm for Signal Segmentation. 
\textit{Entropy}, 20, 1, 55, 1-20, 2018. 

\bibitem{Hubert2018b}  
Hubert, Paulo; Stern, Julio Michael.
Probabilistic Equilibrium: A Review on the Application of MaxEnt to Macroeconomic Models. 
\textit{Springer Proceedings in Mathematics and Statistics}, 239, 187-197, 2018. 

\bibitem{Hubert2019}
Hubert, Paulo; Killick, Rebecca; Chung, Alexandra; Padovese, Linilson R. 
A Bayesian binary algorithm for root mean squared-based acoustic signal segmentation. 
\textit{Journal of the Acoustical Society of America}, 146, 3, 1799-1807, 2019.  

\bibitem{Hulsroj2013}
Hulsroj, Peter. \textit{The Principle of Proportionality}. Springer Netherlands, 2013. 
  
\bibitem{Hutter2005} 
Hutter, Marcus. \textit{Universal Artificial Intelligence: Sequential Decisions based on Algorithmic Probability}. Berlin, Springer, 2005. 
  
\bibitem{Hutter2008} 
Hutter, Marcus. \textit{Predictive Hypothesis Identification}. Canberra ACT, 2008. \    
 
\bibitem{Irony2002}
Irony, Telba Zalkind; Lauretto, Marcelo de Souza; Pereira, Carlos Alberto de Bragan\c{c}a; Stern, Julio Michael. 
A Weibull Wearout Test: Full Bayesian Approach. pp.287-300 in Hayakawa, Y; Irony, T.Z; Xie, M; eds. 
\textit{Systems and Bayesian Reliability}.   
Singapore, World Scientific, 2002. 

\bibitem{Inhasz2010} 
Inhasz, Rafael; Stern, Julio Michael. Emergent Semiotics in Genetic Programming and the Self-Adaptive Semantic Crossover. \textit{Studies in Computational Intelligence}, 314, 381-392, 2010.   

\bibitem{Izbicki2012}
Izbicki, Rafael; Fossaluza, Victor; Hounie, Ana Gabriela; Nakano, Eduardo Yoshio; Pereira, Carlos Alberto de Bragan\c{c}a.  
Testing allele homogeneity: the problem of nested hypotheses. 
\textit{BMC Genetics}, 13, 103, 1-11, 2012. 

\bibitem{Izbicki2015}  
Izbicki, Rafael; Esteves, Luis Gustavo.  Logical consistency in simultaneous statistical test procedures. 
\textit{Logic Journal of the IGPL}, 23, 732-758, 2015. 

\bibitem{Izbicki2022}  
Izbicki, Rafael; Stern, Julio Michael; Stern, Rafael Bassi; Esteves, Lu\'{\i}s Gustavo. 
Logical Coherence in Bayesian
Simultaneous Three-Way Hypothesis Tests. 
Submitted for publication, 2022. 

\bibitem{Jeffreys1961}       
Jeffreys, Harold.  
\textit{Theory of Probability}. Oxford, Clarendon Press. 1961, first ed. 1939.

\bibitem{Johnson2009} 
Johnson, Ria; Chakrabarty, Dalia; O'Sullivan, Ewan;   Raychaudhury, Somak. 
Comparing X-Ray and Dynamical Mass Profiles in the Early-Type Galaxy NGC 4636. 
\textit{The Astrophysical Journal}, 706, 2, 706, 980-994, 2009. 

\bibitem{Kadane2011}   
Kadane, Joseph Born. Principles of Uncertainty.  
NY, Chapman-Hall/CRC, 2011.  

\bibitem{Kadane2016}   
Kadane, Joseph Born. Pragmatics of Uncertainty.  
NY, Chapman-Hall/CRC, 2016.  

\bibitem{Kades1997}
Kades, Eric. 
\textit{The Laws of Complexity and the Complexity of
Laws: The Implications of Computational
Complexity Theory for the Law}. 
College of William and Mary Law School Scholarship Repository, 1997. 

\bibitem{Kaplan1987}          
Kaplan, Stanley; Lin,  James C.  An Improved Condensation Procedure in Discrete Probability Distribution Calculations. 
\textit{Risk Analysis}, 7, 15-19, 1987. 

\bibitem{Kapur1989} 
Kapur, Jagat Narain. 
\textit{Maximum Entropy Models in Science and Engineering}. 
New Delhi, John Wiley, 1989. 

\bibitem{Kapur1992} 
Kapur, Jagat Narain; Kesavan, Hiremagalur Krishnaswamy. \textit{Entropy Optimization and Applications}. Academic Press, 1992 

\bibitem{Karawatzki2005}
Karawatzki, Roman; Leydold, Josef; P\"{o}tzel-berger, Klaus \textit{Automatic Markov Chain Monte Carlo Procedures for Sampling from Multivariate Distributions}. 
\ \texttt{epub.wu.ac.at/id/eprint/1400}  

\bibitem{Kaufmann1977}
Kaufmann, Arnold; Grouchko, Daniel; Cruon, R.  
\textit{Mathematical Models for the Study of the Reliability of Systems}. 
NY, Academic Press, 1977. 



\bibitem{Kelter2020}
Kelter, Riko. Analysis of Bayesian posterior significance and effect size indices for the two-sample t-test to support reproducible medical research. \textit{BMC Medical Research Methodology}, 20, 88, 1-18, 2020 

\bibitem{KelterStern2020}
Kelter, Riko; Stern, Julio Michael. The Full Bayesian Significance Test and the e-value – Foundations, theory and application in the cognitive sciences, 2020. \texttt{arXiv:2006.03334}. 
     
\bibitem{Kelter2021a}
Kelter, Riko. Bayesian Hodges-Lehmann tests for
statistical equivalence in the two-sample
setting: Power analysis, type I error rates and
equivalence boundary selection in
biomedical research. \textit{BMC Medical Research Methodology},  21, 171, 1-26, 2021. 

\bibitem{Kelter2021b}
Kelter, Riko. On the Measure-Theoretic Premises of Bayes Factor and Full Bayesian Significance Tests: A Critical Reevaluation. \textit{Computational Brain \& Behavior}, 2021. 
\texttt{doi:10.1007/s42113-021-00110-5}.  

\bibitem{Klir1988} 
Klir, George J; Folger, Tina. \textit{Fuzzy Sets, Uncertainty, and Information}. Prentice-Hall, 1988. 

\bibitem{Klir1995} 
Klir, George J; Yuan, Bo. \textit{Fuzzy Sets and Fuzzy Logic}. Prentice-Hall, 1995. 

\bibitem{Kokott1998}
Kokott, Juliane. \textit{The Burden of Proof in Comparative
and International Human Rights Law}. The Hague:
Kluwer, 1998. 
   
\bibitem{Kostrzewski2012}
Kostrzewski, Maciej. 
On the Existence of Jumps in Financial Time Series. 
\textit{Acta Physica Polonica B}, 43, 10, 2001-2019, 2012. 

\bibitem{Lauretto2003} 
Lauretto, Marcelo de Souza; Pereira, Carlos Alberto de Bragan\c{c}a; Stern, Julio Michael, Zacks, Shelemyahu.  Full Bayesian Significance Test Applied to Multivariate Normal Structure Models. 
\textit{Brazilian. Journal of Probability and Statistics}, 17, 147-168, 2003. 

\bibitem{Lauretto2005a} 
Lauretto, Marcelo de Souza; Stern, Julio Michael. FBST for Mixture Model Selection. 
\textit{AIP Conference Proceedings}, 803, 121-128, 2005. 

\bibitem{Lauretto2005b} 
Lauretto, Marcelo de Souza; Stern, Julio Michael. Testing Significance in Bayesian Classifiers. \textit{Frontiers in Artificial Intelligence and Applications}, 132, 34-41, 2005. 

\bibitem{Lauretto2007}   
Lauretto, Marcelo de Souza; Faria, Silvio; Pereira, Basilio de Bragan\c{c}a; Pereira, Carlos Alberto de Bragan\c{c}a; Stern, Julio Michael.  The Problem of Separate Hypotheses via Mixtures Models. \textit{AIP Conference Proceedings}, 954, 268-275, 2007.
 
\bibitem{Lauretto2009}
Lauretto, Marcelo de Souza; Nakano, Fabio; Faria, Silvio; Pereira, Carlos A.\,de Bragan\c{c}a; Stern, Julio Michael.
A Straightforward Multiallelic Significance Test for the Hardy-Weinberg Equilibrium Law. 
\textit{Genetics and Molecular Biology}, 32, 3, 619-625, 2009. 


\bibitem{Lawvere1997} 
Lawvere, F. William; Schanuel, Stephen Hoel. 
\textit{Conceptual Mathematics: A First Introduction to Categories}. Cambridge Univ. Press, 1997. 

\bibitem{Lima2014} 
Lima, Adriano R;  Mello, Marcelo; Andreoli,  S\'{e}rgio; Fossaluza, Victor; Ara\'{u}jo, C\'{e}lia de;  Jackowski, Andrea; Bressan, Rodrigo; Mari, Jair. 
The Impact of Healthy Parenting As a Protective Factor for Posttraumatic Stress Disorder in Adulthood: A Case-Control Study. 
\textit{PLOS ONE}, 9, 1, 1-9, e87117, 2014. 
  
\bibitem{Lindley1972}
Lindley, Dennis Victor. \textit{Bayesian Statistics: A Review}. SIAM, 1972. 

\bibitem{Lindley1991}
Lindley, Dennis Victor. Subjective Probability, Decision Analysis and their Legal Consequences. \textit{J.\,R.\,Stat.\,Soc.}, 154, 83-92. 
   
\bibitem{Lipschutz1965} 
Lipschutz, Seymour. \textit{General Topology}. 
NY, McGraw-Hill, 1965. 
 
\bibitem{Loschi2007}
Loschi, Rosangela Helena; Monteiro, Jo\~{a}o V.D; Rocha, Gustavo H.M.A; Mayrink, Vinicius D. 
Testing and Estimating the Non-Disjunction Fraction in Meiosis I using Reference Priors. 
\textit{Biometrical Journal}, 49, 6, 824-839, 2007. 

\bibitem{Loschi2012}
Loschi, Rosangela Helena; Santos, Cristiano C;  Arellano-Valle, Reinaldo B. 
Test procedures based on combination of Bayesian evidences for $H_0$. 
\textit{Brazilian J.\,of Probability and Statistics}, 26, 4, 450-473, 2012. 

\bibitem{Luenberger2008} 
Luenberger, David Gilbert; Ye, Yinyu. \textit{Linear and Nonlinear Programming}. NY, Springer, 2008. 

\bibitem{Luhmann1985}
Luhmann, Niklas.  \textit{A Sociological Theory of Law}. London, Routledge \& Kegan Paul, 1985. 

\bibitem{Luhmann1989}
Luhmann, Niklas.  Ecological Communication. Chicago Univ. Press, 1989. 
  
\bibitem{Madruga2001}  
Madruga, Maria Regina; Esteves, Luis Gustavo;  Wechsler, Sergio. On the Bayesianity of Pereira-Stern Tests. \textit{Test}, 10, 291-299, 2001. 

\bibitem{Madruga2003}   
Madruga, Maria Regina; Pereira, Carlos Alberto de Bragan\c{c}a; Stern, Julio Michael. Bayesian Evidence Test for Precise Hypotheses.
\textit{Journal of Statistical Planning and Inference}, 117, 185-198, 2003.


\bibitem{Mathis2008} 
Mathis, Maria Alice de; Rosario, Maria C.do; Diniz, Juliana Belo; Torres,  Albina R;   Shavitt, Roseli G; Ferr\~{a}o, Ygor A;   Fossaluza, Victor; Pereira, C.A.B; Miguel, Eur\'{\i}pedes Constantino. 
Obsessive-compulsive disorder: Influence of age at onset on comorbidity patterns. 
\textit{European Psychiatry}, 23, 3, 187-194, 2008. 

\bibitem{Maturana1980} 
Maturana, Humberto Romesín; Varela, Francisco  Javier. \textit{Autopoiesis and Cognition. The Realization of the Living}. 
Dordrecht, Reidel, 1980. 

\bibitem{Minoux1986} 
Minoux, Michel. \textit{Mathematical Programming: Theory and Algorithms}. John Willey, 1986. 

\bibitem{Montoya2001}    
Montoya-Delgado, Luis E; Irony, Telba  Zalkind; Pereira, Carlos A. de B; Whittle, Martin R. 
An Unconditional Exact Test for the Hardy-Weinberg Equilibrium Law: Sample-Space Ordering Using the Bayes Factor. 
\textit{Genetics}, 158, 2, 875-883, 2001.  

\bibitem{Maranhao2012}  
Maranh\~{a}o, Viviane de Luca; Lauretto, Marcelo de Souza; Stern, Julio Michael. FBST for Covariance Structures of Generalized Gompertz Models. 
\textit{AIP Conference Proceedings}, 1490, 202-211, 2012.  



\bibitem{Muchmore2016}
Muchmore, Adam I. Uncertainty, Complexity, and Regulatory Design. \textit{Houston Law Review}, 53, 5, 1321-1367, 2016. 
   
\bibitem{Murray2019}  
Murray, Jamine; Webb, Thomas E; Wheatley, Steven. 
\textit{Complexity Theory and Law: Mapping an Emergent Jurisprudence}. London, Routledge, 2019. 
  
\bibitem{Nakano2006} 
Nakano, Fabio; Pereira, Carlos Alberto de Bragan\c{c}a; Stern, Julio Michael; Whittle, Martin R. 
Genuine Bayesian Multiallelic Significance Test for the Hardy-Weinberg Equilibrium Law. 
\textit{Genetics and Molecular Research}, 4, 619-631, 2006. 


\bibitem{Oliveira2018}
Oliveira, Natalia L; Pereira, Carlos A.B; Diniz, Marcio A; Polpo Adriano.
A discussion on significance indices for contingency tables under small sample sizes.
\textit{PLoS ONE}, 13, 8, e0199102, 1-19, 2018.

\bibitem{Patriota2013} 
Patriota, Alexandre Galv\~{a}o.  
A classical measure of evidence for general null hypotheses. 
\textit{Fuzzy Sets and Systems}, 233, 74-88, 2013. 

\bibitem{Patriota2017}
Patriota, Alexandre Galv\~{a}o.  
On Some Assumptions of the Null Hypothesis Statistical Testing. 
\textit{Educational and Psychological Measurement}, 77, 3, 507-528, 2017.  

\bibitem{Pawitan2001}  
Pawitan Yudi. \textit{In All Likelihood: Statistical Modelling and Inference Using Likelihood}.
Oxford University Press, 2001.












\bibitem{PereiraStern1999} 
Pereira, Carlos Alberto de Bragan\c{c}a; Stern, Julio Michael. Evidence and Credibility: Full Bayesian Significance Test for Precise Hypotheses. 
\textit{Entropy}, 1, 99-110, 1999. 

\bibitem{PereiraStern1999b} 
Pereira, Carlos Alberto de Bragan\c{c}a; Stern, Julio Michael. 
A Dynamic Software Certification and Verification Procedure. 
\textit{Annals of the 5th ISAS-SCI}, 2, 426-435, 1999. 

\bibitem{PereiraNS2000} 
Pereira, Carlos Alberto de Bragan\c{c}a; Nakano, Fabio; Stern, Julio Michael. 
Actuarial Analysis via Branching Processes.
\textit{Annals of the 6th ISAS-SCI}, 8, 353-358, 2000. 
   

\bibitem{Pereira2a} 
Pereira, Basilio de Bragan\c{c}a;  Pereira, Carlos Alberto de Bragan\c{c}a.  
A Likelihood Aproach to Diagnostic Tests in Clinical Medicine. 
\textit{RevStat - Statistical Journal}, 3, 1,  77-98, 2005. 

\bibitem{Pereira2008}
Pereira, Carlos Alberto de Bragan\c{c}a; Stern, Julio Michael; Wechsler, Sergio. Can a Significance Test be Genuinely Bayesian? 
\textit{Bayesian Analysis}, 3, 79-100, 2008. 

\bibitem{PereiraStern2008b}
Pereira, Carlos Alberto de Bragan\c{c}a; Stern, Julio Michael. Special Characterizations of Standard Discrete Models. 
\textit{RevStat - Statistical Journal}, 6, 199-230, 2008. 


\bibitem{Pereira2011}
Pereira, Carlos Alberto de Bragan\c{c}a. Full Bayesian Significant Test (FBST). pp.551-554 in Lovric, M. ed. \textit{International Encyclopedia of Statistical Science}, Berlin, Springer, 2011. 
 
\bibitem{Pereira2b} 
Pereira, Basilio de Bragan\c{c}a;  Pereira, Carlos Alberto de Bragan\c{c}a.  
\textit{Model Choice in Nonnested Families}. Berlin, Springer. 2016. 
   
\bibitem{PereiraStern2020}
Pereira, Carlos Alberto de Bragan\c{c}a; Stern, Julio Michael. The e‑value: A fully Bayesian significance measure for precise statistical hypotheses and its research program.  
\textit{S\~{a}o Paulo Journal of Mathematical Sciences}, 2020.  \texttt{doi:10.1007/s40863-020-00171-7}. 

\bibitem{Pflug1996} 
Pflug, Georg Christian. 
\textit{Optimization of Stochastic Models: The Interface Between Simulation and Optimization}. Dordrecht, Kluwer, 1996.    
 
\bibitem{Pigliucci2014} 
Pigliucci, Massimo; Boudry, Maarten. 
Prove it! The Burden of Proof Game in Science vs. Pseudoscience Disputes. 
\textit{Philosophia}, 42, 2, 487-502, 2014. 

\bibitem{Pinto2012}
Pinto, Anna; Ventura, Laura.  
\textit{Approssimazioni Asintotiche di Ordine Elevato per Verifiche d'Ipotesi Bayesiani: Uno Studio per Dati di Sobrevvivenza}. 
Universit\`{a} degli Studi di Padova, Dipartimento di Scienze Statistiche, 2012.    

\bibitem{Piskunov1969}   
Piskunov, Nikilai. \textit{Differential and Integral Calculus}. Moscow, MIR, 1969.  

\bibitem{Priest2001} 
Priest, Graham. \textit{An Introduction to Non-Classical Logics}. Cambridge Univ. Press, 2001. 

\bibitem{Ranzato2018} 
Ranzato, Giulia; Ventura, Laura.    
\textit{Biostatistica Bayesiana con ``Matching Priors''}. 
Universit\`{a} degli Studi di Padova, Dipart.\,di Scienze Statistiche, 2018.   
 
\bibitem{Rasch2000} 
Rasch, William. \textit{Niklas Luhmann's Modernity: Paradoxes of Differentiation}. Stanford Univ. Press, 2000. 
 

\bibitem{Rine1984} 
Rine, David C. \textit{Computer Science and Multiple-Valued Logic}. North-Holland, 1984. 

\bibitem{Rifo2009}
Rifo, Laura Leticia Ramos; Torres, Soledad. Full Bayesian Analysis for a Class of Jump-Diffusion Models.  
\textit{Communications in Statistics - Theory and Methods}, 38, 8, 1262-1271, 2009. 

\bibitem{Rifo2012}
Rifo, Laura Leticia Ramos; Gonz\'{a}lez-L\'{o}pez, Veronica. Full Bayesian Analysis for a Model of Tail Dependence. 
\textit{Communications in Statistics - Theory and Methods}, 41, 22, 4107-4123, 2012.   

\bibitem{Rincon2010}
Rinc\'{o}n, Sonia V. del; Rogers, Jeff; Widschwendter, Martin; Sun, Dahui; Sieburg, Hans B; Spruck, Charles. 
Development and Validation of a Method for Profiling
Post-Translational Modification Activities Using Protein
Microarrays. 
\textit{PLoS ONE}, 5, 6, e11332, 1-11, 2010. 

\bibitem{Ripley1987}
Ripley, Brian David.  Stochastic Simulation. NY, Wiley, 1987. 
       
\bibitem{Rissanen1989} 
Rissanen, Jorma. \textit{Stochastic Complexity in Statistical Inquire}. Singapore, World Scientific, 1989. 
   
\bibitem{Rodrigues2021} 
Rodrigues, Abilio; Bueno-Soler, Juliana; Carnielli, Walter Alexandre. Measuring Evidence: A Probabilistic Approach to an Extension of Belnap-Dunn Logic. 
\textit{Synthese}, 198, 22, 5451-5480, 2021.   

\bibitem{Rodrigues2006}
Rodrigues, Josemar. Full Bayesian Significance Test for Zero-Inflated Distributions. 
\textit{Communications in Statistics - Theory and Methods}, 35,  299-307, 2006.  

\bibitem{Romero1991}
Romero, Carlos. \textit{Handbook of Critical Issues in Goal Programming}. Oxford, Pergamon Press, 1991.

\bibitem{Royall1997}
Royall, Richard. 
\textit{Statistical Evidence: A Likelihood Paradigm}. London, Chapman and Hall. 1997. 

\bibitem{Ruhl2008}
Ruhl, John B. Law's Complexity: A Primer. \textit{Georgia State Univ. Law Review}, 24, 4, 1-28, 2008.  
    
\bibitem{Ruhl2015}
Ruhl, John B; Katz, Daniel Martin. Measuring, Monitoring, and Managing Legal Complexity. \textit{University of Iowa Law Review}, 101, 191-244, 2015. 
   
\bibitem{Ruli2020}    
Ruli, Erlis; Sartori, Nicola; Ventura, Laura. Robust approximate Bayesian inference. 
\textit{Journal of Statistical Planning and Inference}, 205, 10-22, 2020.   
   
\bibitem{Ruli2020b}    
Ruli, Erlis; Ventura, Laura. Can Bayesian, Confidence Distribution and Frequentist Inference Agree?  
\textit{Statistical Methods and Applications}, 2020.  \texttt{doi:10.1007/s10260-020-00520-y}.    

\bibitem{Russell2015} Russell, Stuart. Unifying Logic and Probability. \textit{Communications of the ACM}, 58, 7, 88-97, 2015. 
    
   
\bibitem{Sakamoto1986} 
Sakamoto,Yosiyuki; Ishiguro, Makio; Kitaga-wa, Genshiro. 
\textit{Akaike Information Criterion Statistics}. Dordrecht, Kluwer, 1986.

\bibitem{Salicone2007} 
Salicone, Simona. \textit{Measurement Uncertainty: An Approach via the Mathematical Theory of Evidence}. Springer, 2007. 
 
\bibitem{Santos2020}
Santos, Natalia C.L; Dias, Rosa; Alvesc, Diego; Melo, Brian Ganassin, Maria; Gomes, Luiz; Severid, Willia; Agostinho, Angelo.  
Trophic and limnological changes in highly fragmented rivers predict the decreasing abundance of detritivorous fish. 
\textit{Ecological Indicators}, 110, 105933, 1-8, 2020.  

\bibitem{Schuck1992}
Schuck, Peter H. Legal Complexity: Some Causes, Consequences, and Cures. \textit{Duke Law Journal}, 42, 1, 1-52, 1992. 

\bibitem{Schwartzman2008} 
Schwartzman, Armin; Mascarenhas, Walter F; Taylor, Jonathan E.  
Inference for Eigenvalues and Eigenvectors of Gaussian Symmetric Matrices. \textit{Ann.\,Statist.}, 36, 6, 2886-2919, 2008. 

\bibitem{Searcoid2006}
Searc\'{o}id, M\'{\i}che\'{a}l \'{O}.  
\textit{Metric Spaces}. London,	Springer, 2006.

\bibitem{Seixas2008} 
Seixas, Andre Augusto Anderson; Hounie, Ana G; Fossaluza, Victor; Curi, Mariana; Alvarenga, Pedro G; Mathis,  Maria A; Vallada, Homero; Pauls, David; Pereira, C.A.B; Miguel, Euripedes Constantino.  
Anxiety Disorders and Rheumatic Fever: Is There an Association? 
\textit{CNS Spectrums}, 13, 12, 1039-1046, 2008.  
   
\bibitem{Self1987}  
Self, Steven G; Liang, Kung-Yee. Asymptotic Properties of Maximum Likelihood Estimators and Likelihood Ratio Tests Under Nonstandard Conditions. \textit{J.\,Amer.\,Statist.\,Assoc.}, 82, 398, 605-610, 1987. 
   
\bibitem{Shackle1968} 
Shackle, George Lennox Sharman. 
\textit{Uncertainty in Economics and Other Reflections}. 
London, Cambridge Univ. Press, 1968. 

\bibitem{Shackle1969}
Shackle, George Lennox Sharman.  
\textit{Decision, Order and Time in Human Affairs}. 
London, Cambridge Univ. Press, 1969.    

\bibitem{Shavitt2017} 
Shavitt, Roseli G; Requena, Guaraci; Alonso, Pino; Zai, Gwyneth; Costa, Daniel L.C; Pereira, C.A.B; Ros\'{a}rio, Maria C; Morais, Ivanil; Fontenelle, Leonardo; Cappi, C; Kennedy, James; Menchon, Jose M; Miguel, Eur\'{\i}pides; Richter, Peggy M.A.  Quantifying dimensional severity of obsessive-compulsive disorder for neurobiological research. 
\textit{Progress in Neuro-Psychopharmacology and Biological Psychiatry}, 79, 206-212, 2017. 

\bibitem{Shier1991} 
Shier, Douglas R. \textit{Network Reliability and Algebraic Structures}. Oxford, Clarendon Press, 1991.  

\bibitem{Sichelman2021}
Sichelman, Ted. Quantifying Legal Entropy. 
\textit{Frontiers in Physics}, 21, 9, 665054, 1-14, 2021.  

\bibitem{Silva2015}
Silva, Gustavo Miranda da; Esteves, Luis Gustav; Fossaluza, Victor; Izbicki, Rafael; Wechsler, Sergio. 
A Bayesian Decision Theoretic Approach to Logically-Consistent Hypothesis Testing. 
\textit{Entropy}, 17, 10, 6534-6559, 2015. 

\bibitem{Silva2018} 
Silva, Ivair R. On the correspondence between frequentist and Bayesian tests. 
\textit{Communications in Statistics - Theory and Methods}, 47, 14, 3477-3487, 2018. 

\bibitem{Sikov2019}  
Sikov, Anna; Stern, Julio M. Application of the full Bayesian significance test to model selection under informative sampling. 
\textit{Statistical Papers}, 60, 89-104, 2019. 

\bibitem{Skilling2004}  
Skilling, John. Nested Sampling. \textit{AIP Conference Proceedings}, 735, 395-405, 2004. 

\bibitem{Skilling2006}  
Skilling, John. Nested Sampling for General Bayesian Computation. \textit{Bayesian Analysis}, 1, 4, 833-860, 2006. 

\bibitem{Spall2003} 
Spall, James C. \textit{Introduction to Stochastic Search and Optimization}. Willy, 2003. 

\bibitem{Spektor2019} 
Spektor, Mikhail S; Gluth, Sebastian; Fontanesi, Laura; Rieskamp, J\"{o}rg.  
How Similarity Between Choice Options Affects Decisions From Experience: The Accentuation-of-Differences Model. 
\textit{Psychological Review}, 126, 1, 52-88, 2019. 

\bibitem{Spivak2014}
Spivak, David Isaac. \textit{Category Theory for the Sciences}. MIT Press, 2014. 

\bibitem{Springer1979} 
Springer, Melvin D. \textit{The Algebra of Random Variables}. NY, Wiley, 1979.

\bibitem{Stern1992} 
Stern, Julio Michael. Simulated Annealing with a Temperature Dependent Penalty Function. 
\textit{ORSA Journal on Computing}, 4, 311-319, 1992. 

\bibitem{SternZacks2002}
Stern, Julio Michael; Zacks, Shelemyahu. Testing the Independence of Poisson Variates under the Holgate Bivariate Distribution: The Power of a New Evidence Test. 
\textit{Statistical and Probability Letters}, 60, 313-320, 2002.

\bibitem{Stern2003} 
Stern, Julio Michael. Significance Tests, Belief Calculi, and Burden of Proof in Legal and Scientific Discourse. 
4th Laptec -- Logic Applied to Technology congress,  
Advances in Intelligent Systems and Robotics.   
\textit{Frontiers in Artificial Intelligence and its Applications}, 101, 139-147, 2003. 
   
\bibitem{Stern2004} 
Stern, Julio Michael.  Paraconsistent Sensitivity Analysis for Bayesian Significance Tests.  
\textit{Lecture Notes in Artificial Intelligence}, 3171, 134-143, 2004. 

\bibitem{Stern2007a}
Stern, Julio Michael. Cognitive Constructivism, Eigen-Solutions, and Sharp Statistical Hypotheses. \textit{Cybernetics \& Human Knowing}, 14, 1, 9-36, 2007. 

\bibitem{Stern2007b}  
Stern, Julio Michael. Language and the Self-Reference Paradox. 
\textit{Cybernetics \& Human Knowing}, 14, 4, 71-92, 2007. 

\bibitem{Stern2008a} 
Stern, Julio Michael. Decoupling, Sparsity, 
Randomization, and Objective Bayesian Inference. 
\textit{Cybernetics \& Human Knowing}, 15, 2, 49-68, 2008. 

\bibitem{Stern2008b} 
Stern, Julio Michael. 
\textit{Cognitive Constructivism and the Epistemic Significance of Sharp Statistical Hypotheses in Natural Sciences}.  
Tutorial text for MaxEnt 2008 - The 28th International Workshop on Bayesian Inference and Maximum Entropy Methods in Science and Engineering,   
Borac\'{e}ia, S\~{a}o Paulo, Brazil, July 6-11, 2008.
\texttt{arXiv:1006.5471}.


\bibitem{Stern2011a}
Stern, Julio Michael. Symmetry, Invariance and Ontology in Physics and Statistics. 
\textit{Symmetry}, 3, 3, 611-635, 2011. 

\bibitem{Stern2011b}
Stern, Julio Michael. Constructive Verification, Empirical Induction, and Falibilist Deduction: A Threefold Contrast. 
\textit{Information}, 2, 635-650, 2011. 
  
\bibitem{Stern2014a}
Stern, Julio Michael. Jacob's Ladder and Scientific Ontologies. 
\textit{Cybernetics \& Human Knowing}, 21, 3, 9-43, 2014. 

\bibitem{Stern2014b} 
Stern, Julio Michael; Nakano, Fabio. Optimization Models for Reaction Networks: Information Divergence, Quadratic Programming and Kirchhoff’s Laws. \textit{Axioms}, 3, 109-118, 2014.    

\bibitem{SternPereira2014}
Stern, Julio Michael; Pereira, Carlos Alberto de Bragan\c{c}a.  
Bayesian Epistemic Values: Focus on Surprise, Measure Probability! 
\textit{Logic Journal of the IGPL}, 22, 236-254, 2014. 

\bibitem{Stern2015} 
Stern, Julio Michael. Cognitive-Constructivism, Quine, Dogmas of Empiricism, and M\"{u}nchhausen's Trilemma. Ch.5, pp.55-68, in Polpo, Adriano; Louzada, Francisco; Rifo, Laura L.R;  Stern, Julio M; Lauretto, Marcelo; eds. Interdisciplinary Bayesian Statistics: EBEB 2014, 
\textit{Springer Proceedings in Mathematics and  Statistics}, 118, 2015. 

\bibitem{Stern2017a} 
Stern, Julio Michael. Continuous versions of Haack's Puzzles: Equilibria, Eigen-States and Ontologies. 
\textit{Logic Journal of the IGPL}, 25, 4, 604-631, 2017. 

\bibitem{Stern2017b} 
Stern, Julio Michael. Jacob's Ladder: Logics of Magic, Metaphor and Metaphysics -- Narratives of the Unconscious, the Self, and the Assembly. 
\textit{Sophia}, published at Springer Online First June 7, 2017. \texttt{doi:10.1007/s11841-017-0592-y}. 
 
\bibitem{Stern2018a} 
Stern, Julio Michael. Verstehen (causal/ interpretative understanding), Erkl\"{a}ren (law-governed description/ prediction), and Empirical Legal Studies. 
\textit{Journal of Institutional and Theoretical Economics} (JITE), 174, 105-114, 2018.  

\bibitem{Stern2018b} 
Stern, Julio Michael. Karl Pearson on Causes and  Inverse Probabilities: Renouncing the Bride, Inverted Spinozism and Goodness-of-Fit. 
\textit{South American Journal of Logic}, 4, 1, 219-252, 2018. 

\bibitem{Stern2018c}
Stern, Julio Michael; Izbicki, Rafael; Esteves, Luis Gustavo; Stern, Rafael Bassi. 
Logically-Consistent Hypothesis Testing and the Hexagon of Oppositions.  
\textit{Logic Journal of the IGPL}, 25, 741-757, 2018.  
  
\bibitem{Stern2020} 
Stern, Julio Michael. A Sharper Image: The Quest of Science and Recursive Production of Objective Realities. \textit{Principia}, 24, 2, 255-297, 2020. 
   
\bibitem{Stern2022} 
Stern, Julio Michael. Color-Coded Epistemic Modes in a Jungian Hexagon of Opposition. In Jean-Yves Beziau, Ioannis Vandoulakis (2022), \textit{Studies in Universal Logic}, 20, 303-332, 2022. 
\texttt{doi:10.1007/978-3-030-90823-2\_14} 
   
\bibitem{Stern2023} 
Stern, Julio Michael. Dynamic Oppositional Symmetries for Color, Jungian and Kantian Categories. \textit{Logica Universalis}, 2023.  \texttt{doi:10.1007/s11787-023-00342-y}   

      

\bibitem{Takada2020}
Takada, Hellinton Hatsuo; Ribeiro, Celma de Oliveira; Costa, Oswaldo Luiz do Valle; Stern, Julio Michael.  Gini and Entropy-Based Spread Indexes for Primary Energy Consumption Efficiency and CO2 Emission. \textit{Energies}, 13, 4938, 1-17, 2020. 
   
\bibitem{Teubner1988}
Teubner, Gunther. \textit{Autopoietic Law: A New Approach to Law and Society}. Berlin, Walter de Gruyer, 1988. 
 
\bibitem{Thulin2014} 
Thulin, Mans. Decision-theoretic justifications for Bayesian hypothesis testing using credible sets. \textit{Journal of Statistical Planning and Inference}, 146, 133-138, 2014. 

\bibitem{Ventura2013}  
Ventura, Laura; Ruli, Erlis; Racugno, Walter.  
A note on approximate Bayesian credible sets based on modified loglikelihood ratios. 
\textit{Statistics and Probability Letters}, 83, 11, 2467-2472, 2013.   

\bibitem{Ventura2014} 
Ventura, Laura; Reid, Nancy.   
Approximate Bayesian computation with modified log-likelihood ratios. 
\textit{Metron}, 72, 231-245, 2014. 

\bibitem{Ventura2016} 
Ventura, Laura; Racugno, Walter. Pseudo-Likelihoods for Bayesian Inference. pp.205-220 in DiBattista, T; Moreno, E; Racugno, W; eds. 
\textit{Topics on Methodological and Applied Statistical Inference}. Berlin, Springer, 2016. 

\bibitem{Vieland2019} 
Vieland, Veronica J; Chang, Hasok.  No evidence amalgamation without evidence measurement. 
\textit{Synthese}. 196, 3139-3161, 2019.     

\bibitem{Vikas2016} 
Vikas, K; Rao, Aruna K. 
Full Bayesian Empirical Likelihood Significance Test for
Equality of Medians. 
\textit{InterStat}, 2016, 01, 001, 1-9, 2016.  

\bibitem{Voigt1990} 
Voigt, Hans-Michael; M\"{u}hlenbein, Henz; Schwefel, Hans-Paul. \textit{Evolution and Optimization '89}. 
Berlin, Akademie Verlag, 1990. 

\bibitem{Vosseler2016} 
Vosseler, Alexander; Weber, Enzo. Bayesian analysis of periodic unit roots in the presence of a break. 
\textit{Applied Economics}, 49, 38, 3841-3862, 2016.
  
\bibitem{Wah2008}
Wah, Benjamin W; Chen, Yixin; Wang, Tao. 
Theory and Applications of Simulated Annealing for Nonlinear Constrained Optimization. 
Ch.9, pp.155-186 in Tan, Cher Ming. \textit{Simulated Annealing}, 2008. 
  
\bibitem{Wallace1968}
Wallace, Christopher Stewart; Boulton David Morris  
An Information Measure for Classification. \textit{Computer Journal}, 11, 2, 185-194, 1968. 

\bibitem{Wallace1999}
Wallace, Christopher Stewart; Dowe, David. Minimum Message Length and Kolmogorov Complexity. \textit{Computer Journal}, 42, 4, 270-283, 1999.

    
\bibitem{Wechsler2008} 
Wechsler, Sergio; Pereira, Carlos Alberto de Bragan\c{c}a; Marques, Paulo Cilas. Birnbaum's Theorem Redux. 
\textit{AIP Conf. Proceedings}, 1073, 96-100, 2008. 
 
\bibitem{West1997} 
West, Mike; Harrison, Jeff. \textit{Bayesian Forecasting and Dynamic Models}. Springer, 1997. 

\bibitem{Williams2021}  
Williams, David. \text{Weighing the Odds}. Cambridge University Press, 2001.

\bibitem{Williamson1989}
Williamson, Robert Charles. \textit{Probabilistic Arithmetic}. Ph.D. Thesis, Univ. of Queensland, 1989.

\bibitem{Wittenburg2016}  
Wittenburg, D\"{o}rte; Teuscher, Friedrich;  Klosa, Jan; Reinsch, Norbert. Covariance Between Genotypic Effects and its Use for Genomic Inference in Half-Sib Families.
\textit{G3 - Genes Genomes Genetics}, 6, 9, 2761-2772, 2016.  

\bibitem{Wittgenstein1921} 
Wittgenstein, Ludwig.  
Logisch-Philosophische Abhandlung, \textit{Annalen der Naturphilosophie}, 14, 185-262, 1921. Translated as 
\textit{Tractatus Logico Philosophicus}, NY, Harcourt \& Brace, 1922.   

\bibitem{Zeleny1979} 
Zeleny, Milan. \textit{Autopoiesis, Dissipative Structures and Spontaneous Social Order}. Boulder, CO, AAAS, 1979.  

\bibitem{Zeleny1982} 
Zeleny, Milan. \textit{Multiple Criteria Decision Making}. 
NY, McGraw-Hill, 1982. 

\bibitem{Zellner1971} 
Zellner, Arnold.  
\textit{Introduction to Bayesian Inference in Econometrics}. NY, Wiley, 1971.

\end{thebibliography}
\end{document}